\documentclass [11pt,a4paper,reqno]{amsart}

\usepackage{mathabx,amssymb}

\newcommand{\ds}[1]{\ {#1} \ }
\newcommand{\dss}[1]{\quad {#1} \quad }

\input xy
\xyoption{all}
\def\To{\longrightarrow}

\def\Mdot{\cdot_M}

\def\Ndot{\cdot_N}

\def\iI{{ i \in I }}
\def\jJ{{ j \in J }}

\def\Iff{\ \Leftrightarrow \ }

\def\zt{\zeta}
\def\00{\{ 0\}}

\def\Fate{\operatorname{Fate}}

\def\Eq{\operatorname{Eq}}
\def\Feq{\operatorname{Feq}}

\def\Irr{\operatorname{Irr}}

\def\pipe{\mid}

\def\vrp{\varphi}
\def\m{m}

\def\brv{{\bar v}}

\def\iv{v^{-1}}

\def\Sig{\Sigma}

\def\htU{{\widehat U}}

\def\brU{{\overline U}}

\def\brw{\overline{w}}

\def\ipi{\pi^{-1}}

\def\irho{\rho^{-1}}

\def\sig{\sigma}

\def\sm{\setminus}

\def\onto{\twoheadrightarrow}
\def\Onto{\; -\hskip-5pt\twoheadrightarrow}

\def\bbN{\mathbb N}
\def\bbR{\mathbb R}

\newcommand{\etype}[1]{\renewcommand{\labelenumi}{(#1{enumi})}}

\def\N{\mathbb{N}}
\def\NN{\N_0}
\def\Z{\mathbb{Z}}
\def\Q{\mathbb{Q}}
\def\R{\mathbb{R}}
\def\F{\mathbb{F}}

\def\Eal{E(\al)}

\def\STR{\operatorname{STR}}
\def\STROP{\operatorname{STROP}}
\def\STROPm{\STROP_m}

\def\tT{\mathcal T}

\def\tG{\mathcal G}

\def\MFC{\operatorname{MFC}}

\def\dispace{\setlength{\itemsep}{2pt}}
\def\eroman{\etype{\roman} \dispace}
\def\ealph{\etype{\alph} \dispace}

\def\pSkip{\vskip 1.5mm \noindent}

\def\N{\mathbb N}

\def\mfA{\mathfrak A}

\def\mfz{\mathfrak z}
\def\mfq{\mathfrak q}

\def\mfa{\mathfrak a}

\def\olz{\overline z}
\def\al{\alpha}
\def\bt{\beta}
\def\gm{\gamma}

\def\tlE{\widetilde{E}}

\def\tlU{\widetilde{U}}
\def\tlE{{\widetilde{E}}}

\def\tlal{\widetilde{\al}}

\def\tl0{\widetilde{0}}

\newtheorem{thm}{Theorem} [section]
\newtheorem*{thm*}{Theorem}
\newtheorem{cor}[thm]{Corollary}
\newtheorem{lem}[thm]{Lemma}

\newtheorem{prop}[thm]{Proposition}

\newtheorem*{claim*} {Claim}
\newtheorem*{theorem6.6'} {Theorem 6.6$'$}

\newtheorem*{applc*} {Application}

\newtheorem{acknowledgment*}[thm] {Acknowledgment}

\newtheorem{examp}[thm]{Example}

 \newtheorem*{remark*}{Remark}
\newtheorem{construction}[thm]{Construction}

\newtheorem{schol}[thm]{Scholium}
\newtheorem{notation}[thm]{Notation}

\newtheorem*{notation*} {Notation}
\newtheorem*{comment*} {Comment}
\newtheorem{diagram}[thm]{Diagram}

\theoremstyle{definition}
\newtheorem{defi/}{Definition}

\theoremstyle{definition}
 \newtheorem{defn}[thm]{Definition}
\newtheorem*{defn*} {Definition}
\newtheorem*{example*} {Example}
\newtheorem{examples}[thm]{Examples}
 \newtheorem{rem}[thm]{Remark}
 \newtheorem{rems}[thm]{Remarks}
 \newtheorem{remark}[thm]{Remark}
  
\newtheorem{example}[thm]{Example}

 \renewcommand{\sectionmark}[1]{}

\newcommand{\diag}{\operatorname{diag}}

\newcommand{\Cov}{\operatorname{Cov}}
\newcommand{\Con}{\operatorname{Con}}

\newcommand{\lm}{\lambda}
\newcommand{\Lm}{\Lambda}
\newcommand{\om}{\omega}

 \newcommand{\id}{\operatorname{id}}

\begin{document}

\title[Supertropical Monoids III] {Supertropical Monoids III: \\[1mm] Factorization and splitting covers}
\author[Z. Izhakian]{Zur Izhakian}\address{
Department of Mathematics, Ariel University, 40700, Ariel, Israel.
}
    \email{zzur@g.ariel.ac.il}

    \author[M. Knebusch]{Manfred Knebusch}
\address{Department of Mathematics,
NWF-I Mathematik, Universit\"at Regensburg 93040 Regensburg,
Germany} \email{manfred.knebusch@mathematik.uni-regensburg.de}
%
%
%
%


\date{\today}




\begin{abstract}
The category $\STROP_m$ of supertropical monoids,  whose morphisms are transmissions, has the full--reflective subcategory  $\STROP$ of commutative semirings. In this setup,    quotients are determined directly by equivalence relations, as ideals are not applicable for monoids, leading to a new approach to factorization theory. To this end, tangible factorization into irreducibles is obtained through fiber contractions and their hierarchy. Fiber contractions also provide different quotient structures, associated with covers and types of splitting covers.\end{abstract}

\maketitle
\setcounter{tocdepth}{1}

{\small \tableofcontents}

\baselineskip 14pt

\numberwithin{equation}{section}

\section*{Introduction}
Commutative algebra over  semirings, especially bipotent semirings with respect to monoid valuations (written \m-valuations), plays an important role in classical algebraic geometry \cite{CHWW,Gat,IMS,Mac}. These semirings can be extended to supertropical semirings, over which \m-valuations are refined by supervaluations.
Supervaluations generalize  classical valuations,   having  ordered monoids as their targets, instead of ordered abelian groups ~ \cite[~\S2]{IKR1}, and thereby allow  for  approaching other important valuation targets, such as~ $\bbN$.
 Classical valuations, with~ $\bbR$ as a target, are at the heart of tropical geometry, and have been  motivating our development of a systematic theory of supertropical semirings and  supervaluations,  relayed on a solid algebraic foundation \cite{I1,I2}, \cite{IR1}--\cite{IR6},~ \cite{IKR5}. This paper,  as a sequel of \cite{IKR4,IK},  develops further the theory of supertropical monoids, combining algebraic and categorical viewpoints.
 It utilizes the framework introduced  in ~\cite{IKR1,IKR4}, which can be used as a supporting  glossary, employing only part of the theory developed within these papers. The basics of this framework are described below.

\subsection*{Monoids and semirings in mathematical studies} $ $
\pSkip
 Ordered  monoids,  for example the positive numbers $\Z_{\geq 0}$, $\Q_{\geq 0}$, and $\R_{\geq 0}$ together with the standard summation, play a key role in many classical areas  of study, e.g., combinatorics, number theory, arithmetic geometry, and automata theory, as well as in  more  recent areas, e.g.,  schemes over ~$\F_1$ ~\cite{CC}, monoid schemes \cite{CHWW}, tropical geometry \cite{Mac}, and  blueprints~ \cite{Lor}. Unlike groups, inverse operation does not exist in monoids, and therefore also in semirings,
  so that well-known familiar  methods from commutative algebra are not accessible. Nevertheless,
  every ordered monoid gives rise to a bipotent semiring, which extends naturally to a supertropical semiring (cf. Construction~ \ref{constr:1.2}). Supertropical semirings have a richer algebraic structure that enables recovering of concepts from  commutative algebra, including new methods for analyzing and undersetting monoids. The     present paper continues  to enhance this approach.

Classical commutative algebra over rings lays   the foundations
to algebraic geometry, and to many other areas of mathematics.
To implement  concepts from classical algebra in semiring theory, and thereby in tropical algebra,
one  has to  describe  tropical notions in a formal algebraic language, more amenable to structure theory.
%
%
%
  For example, the  ``ghost copy'' $\tG(R) = A^\nu$ of
   $R = A \cup A^\nu$
 adjoined to the
original max-plus algebra~$A$, in which supertropicality (i.e., $a+a = a^\nu $) replaces  additive idempotence, frames additive doubling. In this setup, the variety of a collection of polynomials is
 the set of elements on which all polynomials are evaluated as a ghost,  namely  $A^\nu$ plays the  role of zero in classical algebra ~\cite{I2}.  This formulates the notion  of ``corner locus'' in
tropical geometry ~\cite{Gat}, and   enables  the recovering of  classical results, for example to utilize the Zariski correspondence  between algebraic sets and geometric  congruences ~\cite{IR6}.

More generally, since  ideals are not applicable over semirings, they cannot be used to build a factor algebra. Aiming to  construct a factor algebra systematically, which also addresses universal properties  in category theory, we combine  equivalence relations (instead of ideals) with the algebraic benefit of the supertropical structure, providing a general suitable framework for commutative algebra over semirings.
In a similar spirit, prime and radical ideals can be replaced by special equivalence relations to develop a systematic theory of  commutative $\nu$-algebra~ \cite{I2}, which generalizes supertropical algebra. These algebraic structures show the advantage of supertropical theory.
They enrich the algebraic language and enable algebraic descriptions of combinatorial properties, and thereby recovering of familiar  mathematical  notions; supertropical  matrices provide  a long list of examples   \cite{IzhakianRowen2008Matrices}--\cite{IzhakianRowen2009Resulatants}.

\subsection*{Factorization in mathematical studies} $ $
\pSkip
A rich theory of factorization into irreducibles has appeared in the literature in the past two centuries. It was initially approached by algebraic number theory, with intensive study of factorization in the setting of Dedekind domains, and later systematically studied in the more general setting of integral domains by  Anderson, Anderson, and Zafrullah in a sequence of papers, starting with \cite{AAZ1} and \cite{AAZ2}. Ever since, a more elaborated systematics theory of factorization has been developed over various algebraic structures, including commutative rings \cite{JM,JMR}, commutative monoids \cite{CCGS,GZ}, commutative semirings \cite{ABLST,CrG}, and even noncommutative structures  \cite{BG,BS}. The foundations  of this theory have  been utilized by the second author in \cite{KK,KZ1,KZ2}. The present paper focuses on particular commutative monoids and semirings.

\subsection*{Supertropical monoids,  supervaluations, and factorizations} $ $
\pSkip Given an  ordered monoid $M = (M,  \cdot \;) $, we define  the bipotent semiring $T(M)$
  with $ \cdot$  as multiplication  and the maximum in the given
ordering as addition. For example, $T(\bbR)$ is the max-plus (tropical) semiring,  obtained by taking the real numbers with the standard ordering and summation. In turn, any bipotent semiring extends to a  supertropical semiring $R$, in which $a + a = \nu(a)$ for any ~$a \in R$, where  $\nu: R \to \tG(R)$  is the   \textbf{ghost map} -- a projection on the    \textbf{ghost ideal}  $\tG(R)$. Equivalently, $a + a = e a$, where $e \in \tG(R) $ is a distinguished idempotent element satisfying $e R = e \tG (R)= \tG (R)$.
The addition of ~$R$,
 which  replaces the idempotent addition  $a+ a = a$ in max-plus algebra, enables  to formulate   combinatorial properties in algebraic manner. In this setup, $\tG(R)$ plays  the role of  zero  in classical algebra. $\tT(R) = R \sm \tG(R)$ is the set of \textbf{tangible} elements in $R$, so that, a  bipotent semiring is  a supertropical semiring without tangible elements.

Tangible elements are the meaningful elements, which essentially  capture  the central algebraic and geometric features within the theory.
Factorization of tangibles
is a delicate  issue, since typically the ghost map $\tT(R) \to \tG(R)$, $x \mapsto ex$, is not injective.
Thus,  the fibers  $x^{-1}(e x)$  need a careful  study,
 including a categorical viewpoint which involves valuations and transmissions, as well as the generalization of supertropical semirings to supertropical monoids. We recall a motivating example
from \cite{IK}.

 \begin{example}\label{exmp:0}
    Consider the  supertropical semiring $R$ with totaly ordered  ghost ideal  $eR = \{e, c, c^2,\dots  \} \ds{\dot\cup} \{ 0 \}$, $c > e$,  all  $c^i$ different, and $\tT(R)$  the free abelian monoid on two generators ~$x,y$. Let $\nu: R \to eR$  be  the ghost map that sends $x^i y^j$ to $c^{i+j} =e x^i y^j$, so that the tangible preimage of~ $c^{n}$ consists of all $x^i y^j$ with $i+j =n$.  Defining the equivalence relation $z \sim z'$ if  $ez = ez'$,  respects the semiring operations. Under this relation   a product ~$z x$ is identified with $z' x$ for every~$x$, and thus factorizations of $z x$ and~    $z' x$ coincide when  quotienting by the   equivalence relation $\sim$.
   Therefore, as factorizations respect quotients, they can be explored in the terms of
factorization of maps that agree with ghost maps, wherein one needs to take special care of the
particular pathologies arising in this setting.
   %
%
\end{example}

Having this  factorization approach,
we point out the links of supertropical theory to classical theory, provided by  \textbf{supervaluations}.
which refine classical valuations \cite{IKR1}--\cite{IKR3} by replacing their target semirings by supertropical semirings.
Let $R$ be a semiring and ~$M$ a bipotent semiring.\footnote{In this paper we  assume that all occurring monoids are commutative.}
A multiplicative monoid homomorphism $v: R \to M $ is an \textbf{\m-valuation}  on  $R$, if it satisfies $ v(x+y) \ds \leq v(x) + v(y)$, cf.  \cite[\S2]{IKR1}. It is a  \textbf{valuation}~$v$, if in addition  $M\sm \{ 0 \}$ is a cancellative multiplicative monoid.
When $R$ is a ring, these valuations are those
 defined by Bourbaki~ \cite{B}. They  provide a mapping of algebraic objects, called \textbf{tropicalization}.
 Not all \m-valuations on rings are valuations, see  \cite[\S 1]{IKR1} for examples. The \m-valuations on semirings are slightly more general than
  those on rings. They  are employed  to obtain the supervaluations.

A \textbf{transmission} $\al:R \to S$  between supertropical semirings is a multiplicative map ~ $\al$ whose restriction $eR \to eS$ to the ghost ideals is a semiring homomorphism (Definition~ \ref{defn1.3}).
Transmissions are more general than semiring homomorphisms \cite{IKR2},  and serve as  morphisms
in the category $\STROP$ of supertropical semirings (Definition \ref{def:ssmr}). They allow to
comprise supervaluations  in  $\STROP$, as their composition $\al \circ v$  with any supervaluation $v$ is again a supervaluation \cite[\S5]{IKR1}. Therefore, $\STROP$ includes bipotent semirings as a subcategory,  providing  a richer algebraic setting  for their study.

%
%

The category $\STROPm$ of \textbf{supertropical monoids}, whose morphisms are transmissions as defined for $\STROP$,  contains $\STROP$ as a full subcategory.  A~supertropical monoid $(U, \cdot \, )$ is a pointed  monoid with 
 a distinguished idempotent element $e = e_U$, for which the subset $eU$ carries a total ordering compatible with multiplication. Then  $eU$ becomes a semiring by defining addition $+: eU \times eU \to eU$  as the maximum
(Definition~ \ref{defn:1.1}).
%
%
This addition extends to whole ~$U$ by using the multiplication~ $\cdot$ and the idempotent element ~$e$. Nevertheless,  distributivity on $(U, \cdot \, , +)$ may fail, but when it does not, $U$ is a supertropical semiring. A morphism $\al: U \to V$ in $\STROPm$  restricts to a homomorphism $eU \to eV$
of bipotent semirings, and obeys the rules of transmissions in \cite[Theorem~5.4]{IKR1}.

A bipotent semiring $M$ can be realized as a supertropical
semiring having only ghost elements, where $e =
1_M$. Then,    the category $\STROP_m /M$ of supertropical monoids over~
$M$ may be viewed as the category of supertropical monoids $U$
having a fixed ghost ideal $eU = M$, whose morphisms
 are  transmissions $\al: U \to V$ with $\al(x) = x$ for all~
$x \in M.$ As in the case of supertropical semirings \cite[\S6]{IKR1},  surjective transmissions over ~$M$ are called \textbf{fiber contractions}.
 When $\al:U \onto V$ is a fiber contraction and $U$ is a
supertropical semiring, then  $V$ is also  a supertropical semiring
 \cite[Theorem 1.6]{IKR4}, and~$\al$ is a semiring homomorphism
\cite[Proposition ~5.10.iii]{IKR1}.

 Any   supertropical monoid $U$ is associated with a supertropical semiring $\htU$, for which there exists a fiber
contraction  $\sig_U: U \to \htU$,   such that every fiber contraction $\al:~U \onto
V$ factors (uniquely) through $\sig_U$, i.e., $\al = \bt \circ \sig_U$, where
 $\bt:\htU \to V$ is  a fiber contraction. Namely, $\STROP /M$ is a full reflective subcategory of $\STROP_m
/M$~\cite[~p9]{F}, \cite[1.813]{FS}.

Universal problems appearing in
$\STROP/ M$  generalize  to $\STROPm /M$ in an obvious way and often lead to easier solutions which are  delivered to $\STROP / M $ by \textbf{reflections} $\sig_U: U \to \htU$.
%
%
%
%
This approach pertains especially  to \m-valuations  $v: R \to M$ on a ring ~$R$, and to  supervaluations $\vrp:R \to S$ with $S$ a supertropical semiring over $e S = M$, as defined in \cite[\S4]{IKR1}.  (They are also defined for $R$ a semiring.)
Such a supervaluation applies to the coefficients of a (Laurent) polynomial $f(\lm) \in R[\lm]$ in a set $\lm$ of $n$
variables, and gives a polynomial $f^\vrp(\lm)$ over $S$.  This view  helps to analyze   supertropical root sets and tangible
components of polynomials $g(\lm) \in F[\lm]$, 
obtained from $f^\vrp(\lm) \in S[\lm]$ by
passing from~ $S$ to various supertropical semifields $F$ \cite[\S5 and~ \S7]{IR1}.

\subsection*{Related structures} $ $
\pSkip
Let $R$ be the supertropical semiring in Example \ref{exmp:0}, with ghost ideal
$M= eR$.
Let $U$ be a supertropical
monoid whose ghost ideal is $M$, and  set $\tT_c(U) := \{z \in \tT (U) \pipe ez = c\}$. Assume that $U$ satisfies the following properties:
\begin{enumerate}\ealph
  \item[(A)] $\tT_c(U)$ is the set of irreducibles of $U$,

  \item[(B)]  $ez = e$ implies $z = 1$ and every tangible element of $U$ factors into irreducible.
\end{enumerate}
In this paper, we  consider a certain supertropical monoid $U$  in which
$eU = M$  to study the factorization of tangible elements in $U$.

Supertropical monoids satisfy properties (A) and (B), which  are therefore relevant to the context of this paper. In modern factorization studies property (B) roughly means that~ $U$ is atomic ~ \cite{G}. The study of atomic commutative monoids and semirings has received intensive attention in recent literature. Therefore, supertropical monoids are far from being an esoteric restriction of atomic monoids, and might be of high interest in future studies \cite{CGGR,CyG}.

\subsection*{Paper outline} $ $
\pSkip
In \cite{IKR4}
 we studied fate distinction, sons, and tyrants in a supertropical monoid $U$. These notions are recalled in \S\ref{sec:8} below, after recalling basics on supertropical semirings and monoids in \S\ref{sec:1}.  We now resume this study with a thrust on explicit examples in the case that $M = e U$ is discrete and totally ordered. Often we assume more specifically that
 \begin{equation}\label{eq:str}
 \begin{array}{ll}
     M = \{e, c, c^2,\dots  \} \ds{\dot\cup} \{ 0 \},  \text{with  $c > 0$,  all  $c^i$ different, and $\tT(R)$ } \\
      \text{is a quotient of  the free abelian monoid $\tT(\mfA)$
    in two letters  ~$u$ and $x$}, \\[2mm]
    U = \mfA / E(\mfa_d), \ \mfa_d \text{ the monoid ideal $\{ u ^{i} x ^j | i+j >d\} $ for some $d>1$.}
\end{array}
     \tag{$*$}
 \end{equation}
 This means that $E(\mfa_d)$ is the fiber contraction of $\mfA$ over $M$ identifying every monomial $u^i x^j$, $i +j >d$, with its ghost companion $e u^i x^j = c^{i+j}$.

 Then we obtain a good hold on fate distinction  in $\tT(U)$,  the sons of all $x \in \tT(U)$, and the
 tyrants and lonely tyrants among these, without resorting to the general theory of equalizers exhibited in \cite[\S6]{IKR4} by use of ``paths''.  Everything can be done here by direct inspection.

 The second main topic of the paper is factorization of tangibles of a supertropical monoid $U$
 without zero divisors into products of irreducible elements.

 When (as before)
  $M = \{e, c, c^2,\dots  \} \ds{\dot\cup} \{ 0 \}$, $c > 0$,  all  $c^i$ different, and $U$  is a supertropical  monoid with $e U = M$, an element $x \in \tT(U)$ is \textbf{irreducible}, if $x \neq 1$ and $x$ is not a product $u v$ of two elements $u,v \in \tT(U) \sm \{ 1 \} $.
  As in the standard theory of monoids, we say that $U$ has \textbf{unique factorization} (written  UF, for short), if any $x \in \tT(U) \sm \{ 1 \} $ can be written as a product of irreducible elements in only one way.

   Given a   set $I$, possibly infinite, let $\mfA = \mfA(I)$ denote the supertropical monoid over ~$M$ with $\tT(\mfA)$ the free abelian semigroup on the letters $x_i $, $i \in I$, and $ex_i = c$ for every $i \in I$. Trivially the $x_i$ are the irreducible elements of $\mfA(I)$ and $\mfA(I)$ has UF. For every
   $\om \in \mfA(I) \sm \{ 1\}$ we obtain in \S\ref{sec:9} a fiber contraction $A_\om$ of $\mfA(I)$ over $M$, in which $\om$ is a lonely tyrant of $A_\om$, the unique one, and the tangibles $z$ of $A_\om$ are monomials in $\tT(\mfA)$ dividing ~$\om$, $\om = z z ^*$ with a unique  $ z ^* \in \tT(A_\om)$.  This paves the way to an understanding of the hierarchy of all $A_\om$ (Proposition \ref{thm:9.3} and Theorems  \ref{thm:9.4}, \ref{thm:9.5}). Each $A_\om$ has UF and fate distinction.

Assuming that $U$ is a \textbf{tangibly finite quotient} of $\mfA = \mfA(I)$ over $M$, i.e., there is a fiber contraction $\pi_U: \mfA \onto U$ over $M$ with $\tT(U)$ finite, we may write $U = \mfA(I)/ E_U$ with $E_U$ an MFCE-relation on $M$. A \textbf{cover} of $U$ is a tangibly finite quotient $W$ of $\mfA(I)$ over $M$ where $E_W \subset E_U$, and so $U \subset W \subset \mfA$, $\pi_U = \pi_{W,U} \circ
\pi _W$, with  a fiber contraction  $\pi_{W,U}: W \onto U$ over $M$. Usually, a given $\zt \in \tT(U)$  does not uniquely  factorize  into irreducibles, i.e., $\ipi_U(\zt)$ is not a singleton in $\tT(\mfA)$.

In the remaining  sections \S\ref{sec:10}--\S\ref{sec:12} we pursue the idea of bringing a given $\zt \in \tT(U)$,  where  $\ipi_U(\zt)$ is not a singleton, nearer to unique factorization into irreducibles by lifting $\zt$ to suitable covers of $U$. We say that a cover $W \supset U$ \textbf{splits $\zt$ totally}, if $\ipi_U(\zt) $ is contained in $\ipi_W(\zt) $, and so $\ipi_U(\zt) = \ipi_W(\zt) $. Clearly, there  is such a unique minimal cover -- the ``\textbf{minimal splitting cover}'' of $U$ for $\zt$.

More comprehensively, there exist partitions $\Pi = (S_j : j \in J)$ of the set $\ipi_U(\zt)$ which admit a cover $W \supset U$ such that $\pi_W(S_j)$ is a singleton $\{ \zt_j\}$ and $\zt_j \neq \zt_k$ for
$j \neq k$. In consequence $\pi_{W,U} (\zt_j) = \zt$ for each $j$. We say that  the cover $W$ is \textbf{adapted} to $\Pi$. These  partitions can (of course) be viewed as equivalence relations on the set $\ipi_U(\zt)$. They are related to fiber contractions over $M$.

We obtain all this more generally  for a family $Z \subset \tT(U)$ instead of a single $\zt \in \tT(U)$, and study \textbf{partial splitting} of $Z$ instead of total splitting. There  we focus on sets $S \subset \ipi_U(Z)$, such that there exists a cover $W \supset U$ with
$\ipi_U(Z) \cap W = S$.

An important issue in this factorization theory is that the setup is very rigid. From $\mfA(I) \supset M$ and the equivalence $E_U$, giving a fiber contraction over $M$, we obtain the entire data without further choices. No auxiliary sets or functions are needed.

\section{Preliminaries}\label{sec:1}

We recall the basics of our  underlying structures.

\begin{defn}\label{def:ssmr}
A \textbf{supertropical semiring} is  a semiring $R = (R,+, \cdot \; )$ where $e:=1+1$ is  additively idempotent
(i.e., $e+e = e$) such that,   for
all~ $a,b\in R$,  $a+b\in\{a,b\}$ whenever $ea \ne eb$ and  $a+b=ea$ otherwise. This implies
 $ e a = 0 \ds \Rightarrow a = 0.$
The \textbf{ghost map}
$\nu: a\mapsto ea$ is a projection on the \textbf{ghost ideal} $eR$ of $R$, which is a bipotent semiring (with unit element ~ $e$),  i.e., $a+b \in \{ a,b\}$, for any $a,b\in eR$.
The semiring $eR$ is totally ordered  by the rule
\begin{equation}\label{eq:1.1}
 a\le b \dss\Leftrightarrow a+b=b.
\end{equation}
$\tT(R):=R\sm(eR)$ denotes  the set of \textbf{tangible} elements of $R$, and $\tG(R) := eR \sm \{ 0 \}$ denotes the set of \textbf{nonzero ghost elements}.
 The
zero $0 = e0$ is regarded mainly as a ghost.
\end{defn}
Taking  a wider viewpoint, we have the following structure.
\begin{defn}\label{defn:1.1}
 A \textbf{supertropical monoid} $U = (U, \cdot \, )$ is
an abelian monoid   with  an
absorbing element $0 := 0_U$, i.e., $0 \cdot x = 0$ for every $x
\in U$, and a distinguished idempotent $e := e_U$ such
that
$$ \forall x \in U: \quad ex =0 \ \Rightarrow \ x = 0.$$
The \textbf{ghost ideal} $M := eU$ of $U$ is totally  ordered by rule \eqref{eq:1.1}, and this ordering  is compatible with multiplication \cite[Definition 1.1]{IKR4}.
The \textbf{ghost map} of \ ~$U$, $\nu_U: U \to M$, $x \mapsto ex$, is a monoid homomorphism. The tangible set  $\tT := \tT(U)$ and the ghost set  $\tG := \tG(U)$ are defined  as in Definition~ \ref{def:ssmr}.
$U$ is called \textbf{unfolded}, if the set
$\tT(U)_0 := \tT(U) \cup ~\{0 \}$ is closed under multiplication.
\end{defn}


An \textbf{MFCE-relation} on a supertropical monoid $U$, with  ghost ideal  $M:= eU$, is  an
equivalence relation $E$ on~$U$, which is multiplicative, and is
fiber conserving, i.e., $x \sim_E y \Rightarrow ex = ey$. Then the  obvious identification  $M/E = M$ holds, and $E$ is a
TE-relation \cite[Definition~ 1.7]{IKR4}.

An element $z \in \tT(U)$ is called a \textbf{son of $x \in \tT(U)$ over $c \in e U$}, if $ez = c$ and there exists $u \in \tT(U)$ such that $z = xu$. Taking the patriarchal  viewpoint,  $z$~ is regarded as the outcome~ $xu$ of pairing $x$ to some $u \in \tT(U)$ such that $exu = c$, without specifying explicitly the ``mother'' element $u$ of $z$.
An element $x \in \tT(U)$ is called a \textbf{tyrant} (or said to be tyrannic), if it has at most one son over any $c \in eU$.
It is said to be \textbf{isolated}, if it has no son~ $z \neq x $ with $ex = ez$. The latter property is weaker than tyrannic, but in our opinion deserves similar attention.

Given a  tyrant $x$ in a supertropical monoid $U$,  the set $\mfa = eU \cup xU$ is  an ideal of~ $U$ containing $M =e U$. This ideal is small in the sense, that each fiber of the restricted ghost  map $\nu_U | \mfa$ has at most two elements, $\mfa_c = \{c,z \}$ if $x$ has a son $z$ over $c$, otherwise $\mfa_c = \{c\}$. Conversely, if  an ideal $\mfa \supset M$ is given with $|\mfa_c| \leq 2$ for every $c \in M$, then $\nu_U: U \to M$ has a section $s: M \to U$ for which $\mfa = M \cup s(M)$. By  \cite[\S5]{IK} it follows that the MFCE-relations, for which every class contains at most one tangible element, are the compressions of these ideals $M \cup s(M)$, cf. \cite[~Definition~ 2.5]{IKR4}.

In an unfolded  monoid $U$,  the set  $N:= \tT(U)_0$ is a multiplicative monoid with absorbing  element~$0$, while  $M := eU$ is a
totally ordered monoid with absorbing  element $0$. The
restriction
$$ \rho := \nu_U |N : N \ds \To M$$
is a monoid homomorphism with $\irho (0) = \{ 0 \}$. Since  $e_U = 1_M = \rho(1_N)$, the supertropical
monoid $U$ is  completely determined by the triple $(N,M, \rho)$.
This allows for constructing   all unfolded supertropical
monoids up to isomorphism.
\begin{construction}[{\cite[Construction 1.3]{IKR4}}] \label{constr:1.2} Let $M$ be a
totally ordered monoid with absorbing element $0_M \leq x$ for
all $x \in M$, i.e., $M$ is can be realized as a bipotent semiring. Let $N$ be a
commutative monoid with  absorbing element  $0_N$ and  a
multiplicative map $ \rho : N \to M$ with $\rho (1_N) =
1_M$, $\irho(0_M) = \{ 0 _N\} $. The set  $U$ is defined to be the disjoint
union of $M \sm \{ 0_M\}$ and  $N \sm \{ 0_N\}$ together with  a new element
$0$, identified as $0 = 0_M = 0_N $, so that
$$U = M \cup N, \quad \text{where } M \cap N = \{ 0 \}.$$
The multiplication on $U$ is given by the rules, in obvious
notation,
$$  x \cdot y = \left\{
\begin{array}{lll}
  x \Ndot y & \text{if} &  x,y \in N, \\[1mm]
  \rho(x) \Mdot y & \text{if} &  x \in N, y \in M, \\[1mm]
  x \Mdot \rho(y) & \text{if} &  x \in M, y \in N, \\[1mm]
  x \Mdot y & \text{if} &  x,y \in M.  \\
\end{array}
\right.
$$
The monoid
 $(U, \cdot \,)$ is then  a (commutative)  unfolded
supertropical monoid
 with $1_U = 1_N$ and absorbing element $0$. Set $e := 1_M$. Then $eU = M$ and $\rho(x) = ex$ for $x \in M $, and for any $x \in U$,  $ex =
0$ iff $x=0$, since $\irho(0) = \{ 0\}$. Thus
$(U, \cdot \, , e)$, together with the ordering on $M = eU$,
is an unfolded supertropical monoid. We denote the
supertropical monoid $U$ by $\STR(N,M, \rho)$.
\end{construction}

This construction generalizes the construction of supertropical
domains \cite[Construction 3.16]{IKR1}  \footnote{The present notation $\STR(N,M, \rho)$ differs slightly from the
notation $\STR(\tT, \tG, v)$ in \cite[Construction ~3.16]{IKR1}, but this causes no confusion
regarding the ambient context.}. There,
to obtain all supertropical predomains up to
isomorphism, it was assumed that $N \sm \{ 0 \}$ and $M \sm \{ 0 \}$ were  closed under
multiplication, and that the monoid $M \sm\{ 0 \}$ is
cancellative. Omitting only the cancellation hypothesis would
give us a class of supertropical monoids not broad enough for our
work below.

Given an \m-valuation $v: R \to M$
with support $\iv(0) = \mfq$, the supertropical semiring~ $U^0(v)$
appearing  in \cite[Theorem 7.4]{IKR4} may be viewed as an instance
of Construction ~\ref{constr:1.2}, as follows.
\begin{example}\label{exm:1.4} Let $E$ denote the equivalence
relation on $R$ having the equivalence classes $[0]_E = \mfq$ and $[x]_E
= \{ x \} $ for $x \in R\sm \mfq$. $E$ is multiplicative, and  hence
 $R / E$ is  a monoid with absorbing element $[0]_E = 0$,
 identified  with the subset $(R \sm \mfq) \cup \{ 0 \}$ in
the obvious~way.  The map $v: R\to M$ induces a monoid
homomorphism $\brv: R/E \to M$ given by $\brv(x) = v(x)$ for $x
\in R \sm \mfq$, $\brv(0) = 0 $. Thus  $\brv^{-1} (0) = \{ 0 \}$
and
$$U^0(v) = \STR(R / E, M , \brv).$$
\end{example}

Many other  constructions and definitions in the setup of \cite{IKR1} and \cite{IKR2} hold when replacing  supertropical semirings by supertropical
monoids, in particular the following one (cf. \cite[\S5]{IKR1})).
\begin{defn}\label{defn1.3}
Let $U$ and $V$ be supertropical monoids.

\begin{enumerate}\ealph
  \item
  A map $\al: U \to V$  is called a \textbf{transmission}, if the following holds:
(i) $\alpha(0)=0$, (ii) $\alpha(1)=1$, (iii) $\forall x,y\in U:
\alpha(xy)=\alpha(x)\alpha(y),$  (iv) $\alpha(e_U)=e_V,$
(v) $\forall x,y\in eU:
 x \leq y  \Rightarrow \alpha(x)  \leq \alpha(y)$.
%
  Namely, $\al$ is a monoid homomorphism sending $0$ to~ $0$
and $e$ to $e$, which restricts to a  homomorphism $\gm = \al^\nu: eU \to
eV$ of bipotent semirings, called the \textbf{ghost part} of~$\al$. We  say that~$\al$ \textbf{covers}
$\gm$.

\item
 The \textbf{zero kernel} of $\al$ is the set
    $ \mfz_\al := \{ x \in U \ | \ \al(x) =0\}.$
    The \textbf{ghost kernel} of $\al$ is the set
    $ \mfa_\al := \{ x \in U \ | \ \al(x) \in eV \}.$

\item
A transmission $\al: U \to V$ is
\textbf{tangible},   if
$\al(\tT(U)) \ds  \subset \tT(V)_0,$
and  is \textbf{strictly tangible}, if
$\al(\tT(U))  \ds \subset \tT(V).$


\end{enumerate}
\end{defn}


The category of commutative semirings, whose morphisms are transmissions, is denoted by  $\STROP$. It is a full and reflective  subcategory of the category of supertropical monoids, denoted by  $\STROP_m$.

For unfolded
supertropical monoids we have an explicit  description of tangible  transmissions.

\begin{prop}[{\cite[Proposition  1.4]{IKR4}}]\label{prop:1.3} Assume that $U' = \STR(N', M', \rho')$
and $U = \STR(N, M, \rho)$ are unfolded supertropical monoids.
\begin{enumerate} \eroman
    \item If $\lm: N' \to N$ is a monoid homomorphism with $\lm(0) =
    0$,  $\mu: M' \to M$ is a semiring homomorphism, and  $\rho' \lm = \mu
    \rho$, then the well-defined map
    $$ \STR(\lm,\mu) : U' = N' \cup M' \dss \To  U = N \cup M,$$
    sending  $x' \in N'$ to $\lm(x')$ and $y' \in M'$ to
    $\mu(y')$,
    is a tangible transmission.

    \item All tangible transmissions from
    $\al : U' \to U$ are obtained in this way.
\end{enumerate}
\end{prop}

%
%

%
%
%

We exhibit some transmissions between supertropical monoids of simple nature \cite{IK}.
\begin{defn}\label{defn2.1}
A
surjective transmission $\al: U \to V$ is called  a \textbf{fiber
contraction}, if the ghost part $\al^\nu = \gm: M \to N$ is an
isomorphism. A  transmission $\al$ is called a \textbf{fiber contraction
over} $M$, if~ $N = M$ and $\gm  = \id_M$.
\end{defn}

Note  that $\al$ is a fiber contraction iff the equivalence
relation $\Eal$ determined by $\al$  is an MFCE-relation,  $\al = \rho \circ
\pi_E$ with $E$ an MFCE-relation,   $\pi_E: U \to U/E$ the canonical projection associated to $E$,  and $\rho$ an isomorphism. Then
$\al$ is a fiber contraction over $M$ iff $M=N$ and $\rho$ is an
isomorphism over~ $M$.

 \begin{defn}\label{defn2.2}
A surjective transmission $\al: U \to V$ is called an \textbf{ideal
compression}, if~ $\al$ is a fiber contraction over $M$ which maps
$U \setminus \mfa_\al$ bijectively  to $V \setminus N = \tT(V)$.
\end{defn}

This means that $\al = \rho \circ \pi_{E(U, \mfa)}$, where  $\mfa$ is an
ideal of $U$ containing the ghost ideal $M$ and $\rho$ an isomorphism from $\brU
:= U/ E(U,\mfa) $ onto $V$ over $M$. Thus  $\mfa = \mfa_\al$.

%
%

If $E$ is an MFCE-relation on $U$, then
$\pi_E: U \to U/E$ is tangible iff ~$E$ is \textbf{ghost
separating} \cite[Definition 6.19]{IKR2}, in other terms,
iff~ $E$ is finer  than the equivalence relation $E_t := E_{t,U}$
on $U$ which has the equivalence classes~$\{ a\in \tT(U) \pipe ex
= a \}$,  $a \in M \sm \{ 0\}$, and the one-point equivalence
classes~$\{y\}$, $y \in M$ \cite[Example 6.4.v]{IKR1}.

\begin{defn}\label{defn2.4}
The MFCE-relations $E$ on $U$ with $E \subset E_t$ are  called
\textbf{tangible MFCE-relations.}
\end{defn}

In this terminology,  the tangible fiber contractions $\al: U \to
V$ over $M := eU$ are the products
$ \al = \rho \circ \pi_T$
with $T$ a tangible MFCE-relation on $U$ and $\rho$ an isomorphism
over ~$M$.

\begin{defn}\label{defn2.5}
A transmission $\al: U \to V$ is called a \textbf{ghost
contraction}, if $\al^\nu$ is a homomorphism from $M$ onto $N$,
and if $\al$ maps $U \sm (M \cup \mfz_\al)$ bijectively onto
$\tT(V) = V \sm N$.
\end{defn}

That is,
$ \al = \rho \circ \pi_{F(U,\gm)}$
with $\gm: M \to N$ a surjective homomorphism, namely $\gm =
\al^\nu$, and $\rho$ an isomorphism over $N$ from $U / F(U,\gm)$
to $V$. Thus $\al$ is a ghost  contraction iff $\al$ is a
surjective pushout transmission in $\STROP_m$. \{The equivalence
relation $F(U,\gm)$ had been introduced in \cite[Theorem 1.13]{IKR4}.\}

%

\section{Supertropical monoids with fate distinction} \label{sec:8}

We begin with a central notion in this paper.

\begin{defn}\label{def:8.1} Let $U$ be a supertropical monoid.
\begin{enumerate}\ealph
  \item We say that two tangible elements $x_1,x_2 \in \tT(U)$ have the \textbf{same fate}, if $ex_1 = ex_2$ and for every $u \in \tT(U)$ the elements $x_1 u $ and $x_2 u$ are of the \textbf{same kind}, i.e., $x_1 u $ and $x_2 u$ are either elements of $\tT(U)$ or of $\tG(U)$, or they are both zero.
      Otherwise, we say that $x_1$ and $x_2$ have \textbf{different fate}.
  \item We say that $U$ has \textbf{fate distinction}, if different tangibles $x_1, x_2$ of $U$ have different fates.

\end{enumerate}
\end{defn}
If $U$ has fate distinction, then suitable quotients of $U$ have this property as well.

\begin{lem}\label{lem:8.2}
  Assume that $U$ has fate distinction, and that $E$ is an MFCE-relation on $U$, for which $\pi_E: U \to U/E$ is a \textbf{strictly tangible fiber contraction}, i.e., $\pi_E$ is a fiber contraction over $M$ with
  $\pi_E(\tT(U)) \subset \tT(U/E)$, cf. \cite[Definitions 2.1 and 2.3]{IKR4}. Then $U/E$ has fate distinction.

\end{lem}
\begin{proof} We identify $M= eU = e(U/E)$, and so have $e[x]_E = [ex]_E = ex$ for every $x \in U$. Given $x_1, x_2 \in \tT(U)$ with $[x_1]_E$ and $[x_2]_E$ different, but $e[x_1]_E = e[x_2]_E$, we conclude that $x_1 \neq x_2$, $x_1 \in \tT(E)$, $x_2 \in \tT(E)$, and $ex_1 = e x_2$. Thus there exists $u \in \tT(U)$ for which   $x_1 u$ and $x_2u$  are of different kind. Consequently,
$[x_1]_E [u]_E = [x_1u]_E$ and $[x_2]_E [u]_E = [x_2u]_E$ are of different kind. This proves that $[x_1]_E$ and   $[x_2]_E$ have different fate.
\end{proof}
We now look for  a ``universal''  MFCE-relation $E$ on $U$ over $M$ such that
$U/E$ has fate distinction.
\begin{defn}\label{def:8.3} We define a binary relation $F$ on $U$ as follows. Let $x_1,x_2 \in U$. If $x_1$ and $x_2$ are tangible, then $x_1 \sim_F x_2$ iff $x_1$ and $x_2$ have the same fate. Otherwise  $x_1 \sim_F x_2$ iff $x_1 = x_2$. We call $F$ the \textbf{equal fate relation} on $U$.
\end{defn}

\begin{thm} The equal  fate relation $F$ is an MFCE-relation on $U$ over $M = eU$. If $U \neq M$, then $U/F  \neq M$ and $\pi_F: U \to U/F$ is a strictly tangible fiber contraction over~ $M$. Every other such contraction $F'$ is coarser than $F$, i.e., $F \subset F'$. (If $U = M$, then of course $\pi_F = \id_M$.)

\end{thm}

\begin{proof} a) Given $x_1, x_2, x_3 \in U$ with $x_1 \sim_F x_2$, $x_2 \sim_F x_3$, we verify that $x_1 \sim_F x_3$. This is evident in the case that $x_1 = x_2$, $x_2 = x_3$, or $x_1 = x_3$. Otherwise $x_1, x_2, x_3$ are tangible. Then we conclude from $x_1 \sim_F x_2$, that $x_1$ and $x_2$ have the same fate, and from $x_2 \sim_F x_3$ that  $x_2$ and $x_3$ have the same fate. Thus $x_1$ and $x_3$ have the same fate, whence $x_1 \sim_F x_3$.

This shows  that $F$ is an equivalence relation on the set $U$. Furthermore,  $x_1 \sim_F x_2$ implies  $ex_1 = ex_2$ for all $x_1, x_2 \in U$, directly by Definition \ref{def:8.3}.

\pSkip
b) Given $x_1, x_2, x_3 \in U$ with $x_1 \sim_F x_2$, we verify that $x_1y \sim_F x_2y$.
We have $e x_1 = e x_2$, whence $ex_1 y = e x_2 y$.  If $x_1 \notin \tT(U)$ or $x_2 \notin \tT(U)$, then $x_1 = x_2$, and so $x_1 y = x_2 y$. Assume  that $x_1 \in \tT(U)$ and $x_2 \in \tT(U)$. Then $x_1$ and $x_2$ have the same fate. If $y \in M$, then $x_1y = e x_1 y = e x_2 y = x_2 y$. In all these cases the claim is evident.

We are left with the case that $y \in \tT(U)$. Then we conclude from $x_1 \sim_F x_2$ that $x_1 y $ and $ x_2 y$
are of the same kind. If these products are in $\tG(U) \cup \{ 0\}$, then $x_1y = e x_1 y = e x_2 y = x_2 y$.
Otherwise $x_1 y $ and $x_2 y$ are in $\tT(U)$. Let $u \in \tT(U)$ be given. If $yu \in M$, then $x_1 y u = ex_1 yu = ex_2 yu  = x_2 yu  $.
If $yu \in \tT(U)$, then $x_1 yu$ and $x_2 yu$ are of the same kind, since $x_1 \sim_F x_2$. This proves that $x_1 y $ and $x_2 y $ have the same fate, whence $x_1y \sim _F x_2 y$.

This shows that the relation~ $F$ is a multiplicative, whence is an MFCE-relation over $M$, and we infer from Definition \ref{def:8.3} that $\pi_F$ is a strictly tangible fiber contraction over~ $M$.

\pSkip
c) Let $F'$  be an MFCE-relation on $U$ over $M$ for which $\pi_{F'}: U \to U/F'$ is a strictly tangible fiber contraction such that  $U/F'$ has fate distinction. Let $x_1, x_2 \in U$ with $x_1 \not \sim_{F'} x_2$.  Then $[x_1]_{F'}$ and
$[x_2]_{F'}$ have different fate, and so $x_1$ and $ x_2$ also have different fate. This proves that $x_1 \not\sim_{F'} x_2 \Rightarrow x_1 \not\sim_{F} x_2  $, whence $x_1 \sim_{F} x_2 \Rightarrow x_1 \sim_{F'} x_2 $, i.e., $F \subset F'$.
%
\end{proof}


We construct elementary examples of supertropical monoids with fate distinction, starting with the supertropical monoid $\mfA$, whose tangible part $\tT(\mfA)$ is the free abelian monoid generated by two letters $u$ and $x$:
\begin{equation}\label{eq:8.a}
 \tT(\mfA) = \{ u^i x^j \pipe i ,j \in \N_0 \}, \qquad \N_0:= \N \cup \{ 0 \},
\end{equation}
consisting of all different words $u^i x^j$, and whose  ghost submonoid is $e \mfA = \tG(\mfA) \dot \cup \{ 0\}$, where
\begin{equation}\label{eq:8.b}
\tG (\mfA) = \{ e, c, c^2, \dots \} \end{equation} is the cyclic ordered semigroup with generator $c > e$, and  $eu = e$, $ex = c$.

\begin{examples}\label{examp:8.5}
  Choose a number $d \in \N$, and  consider the  monomial ideal
  $$ \mfa = \mfa_d := \{ u^i x^j\pipe i + j > d \} \cup \{ c^k \pipe k >d \} \cup \{ 0 \}  $$
  of $\mfA$.
  Let $E := E(\mfa_d)$ denote the equivalence relation on $\mfA$ defined by
  $$ x \sim_E y \ds \Leftrightarrow \text{ either } x=y, \text{ or } ex > c^d, ey > c^d. $$
(The ideal $\mfa_d$ is saturated in the terminology of \cite[Definition 1.15]{IKR4}.) It is clear by
\cite[~\S1]{IKR4}, but can also be verified directly,  that $E$ is an MFCE-relation on $\mfA$ over
$e \mfA$, and so gives a supertropical monoid
$$U_d := \mfA / E(\mfa_d).$$
Furthermore, $\pi_E: \mfA \twoheadrightarrow U_d$ has the zero kernel $\mfa_d$, in fact this is the universal MFCE-relation with this property.

Writing abusively $z$ for a nonzero class $[z]_E$, we obtain
$$ U_d = \{ u^i x^j\pipe 0 \leq i +  j \leq  d \} \cup \{ e, c, \dots, c^d \} \cup \{ 0 \}  $$
with $eu = e$, $ex = c,$ where all elements written down are different. It is easily  verified that the elements $u^i x^{j}$ with $0 \leq i + j \leq d$ have different fates, and so $U_d$ has fate distinction.
To visually demonstrate this, we depict $U_d$ for $d=3$ by the following diagram. (Starting  with the one-element column $j=d$,  proceed to $j = d-1$, etc.)

$$
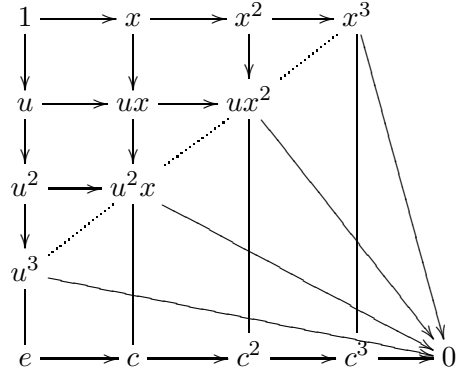
\xymatrix@R=1.5em@C=1.8em{  1  \ar@{>}[r] \ar@{>}[d]&  x \ar@{>}[r]  \ar@{>}[d] & x^2 \ar@{>}[r] \ar@{>}[d] & x^3
\ar @{.}[dl]  \ar@{-}[dddd]  \ar@{>}[ddddr] \\
 u  \ar@{>}[r] \ar@{>}[d] &  ux \ar@{>}[r] \ar@{>}[d] & u x^2 \ar@{.}[ld] \ar@{-}[ddd] \ar@{>}[dddrr]  &
 \\
  u^2  \ar@{>}[r]\ar@{>}[d] &  u^2 x \ar@{.}[dl]  \ar@{-}[dd] \ar@{>}[ddrrr]  &  &
 \\
 u^3 \ar@{-}[d]\ar@{>}[drrrr]  & & & &
 \\
e  \ar @{>}[r]  & c  \ar @{>}[r] & c^2  \ar @{>}[r] & c^3 \ar @{>}[r] & 0  }
$$
\underline{Legend:}
\begin{itemize}
  \item[--] Vertical arrows: multiplication by $u$.
  \item[--] Horizontal arrows: multiplication by $x$.
  \item[--] Oblique  arrows: multiplication by $x$ and by $u$.
  \item[--] Vertical edges: multiplication by $e$.
  \item[--] Dotted line: tangible element without proper sons.

\end{itemize}
\end{examples}

So far we only introduced the terms ``equal fate'' and ``different fate'' for two tangibles in a supertropical monoid  $U$ (Definition \ref{def:8.1}), but we have not yet defined the concept of fate of a tangible element in $U$. This is remedied next.

\begin{defn}\label{def:8.6}
Let $U$ be a supertropical monoid,  $M = eU$. For any $x \in \tT(U)$ we define the function
$$ F_x : \tT(U) \To \{ 1,e,0 \} $$
by
\begin{equation}\label{eq:8.1}
  F_x(u) = \left\{
             \begin{array}{ll}
               1 & \hbox{if } \  xu \in \tT(U), \\
               e & \hbox{if } \ xu \in \tG(U), \\
               0 & \hbox{if } \ xu = 0.
             \end{array}
           \right.
\end{equation}
The \textbf{fate} of $x$ is defined as
\begin{equation}\label{eq:8.2}
\Fate(x) := (ex, F_x).
\end{equation}
We call $ex$ the \textbf{birthday} of $x$.
\end{defn}

The function $F_x$ stores the outcome of the products of $x$ with  elements $u ~\in \tT(U)$, whence the word ``fate''.

\begin{remark} The target $\{ 1,e,0\}$ of the function $F_x$ is a supertropical monoid. More precisely it is the initial object of the category $\STROP_m$ of supertropical monoids, whose morphisms are transmissions. Therefore it is tempting to interpret $F_x$ in terms of this category. Using ~\eqref{eq:8.1}, the function $F_x$ naturally extends to a function $\Phi_x: U \to \{ 1,e,0\} $ on the entire monoid~ $U$. But these functions are  by no mean transmissions, since multiplicativity fails violently. This makes the study of fates intricate from a categorial viewpoint.
\end{remark}

Next we exhibit a second class of examples of supertropical monoids with fate distinction, which in contrast to Examples \ref{examp:8.5} do not have zero divisors. We start with the submonoid $\mfA'$ of the monoid $\mfA$ from above (cf. \eqref{eq:8.a} and \eqref{eq:8.b}), whose underlying set is
$$
\begin{array}{ll}
  \mfA' & = \mfA \sm \{ u^i \pipe i \in \N \} \\[1mm]
   & = \{ 1 \} \cup \{ u^i x^j \pipe i \geq 0, j > 0\} \cup e \mfA \ .
\end{array}
$$
The submonoid
$\mfA'$ is again supertropical, and $e \mfA'$ coincides with $M : = e \mfA$.
We pick a pair $(i,j) \in \NN \times \N$, and introduce the principal ideal
$$ \mfa_{i,j} := u^i x^j \mfA' $$
and the equalizer
$$ E := E_{i,j}= \Eq(u^i x^j, c^j)$$
 in $\mfA' $, which identifies $u^i x^j$ with $c^j$ and so every element $z$ of $\mfa_{i,j}$ with its image $ez \in M$.
 Every  class $[u^k x ^\ell]_E$ with $u^k x^\ell \notin \mfa_{i,j}$ consists only  of the element $u^k x^\ell$.
  This  follows from our theory, since $E$ is an MFCE-relation over $M$, but can also be  verified directly.   Such a class $[u^k x ^\ell]_E$ is abusively denoted as $u^k x^\ell$.
A  class $[u^k x ^\ell]_E$ with $u^k x^\ell \in \mfa_{i,j}$ is usually denoted by $c^\ell \in M$.

Let $$ U'_{i,j} = \mfA/E_{i,j}.$$
 It turns out that $x, ux, u^2x, \dots, u^i x, u^{i+1}x$ have different fates in $U'_{i,j}$, while all elements
 $u^{i+1}x$, $u^{i+2}x, \dots$ have the same fate. Instead of proving this formally, we draw the following diagram of tangibles of $U'_{i,j}$ in the case $(i,j) = (2,3)$, which is used to introduce some arguments concerning the fates.
\begin{diagram}\label{diag:8.7}
 $$
\xymatrix@R=1.5em@C=1.8em{  1  \ar@{>}[r] \ar@{>}[dr]&  x \ar@{>}[r]  \ar@{>}[rd] \ar@{-}[d] & x^2 \ar@{>}[r] \ar@{>}[rd] \ar@{-}[d] & x^3 \ar@{>}[r] \ar@{>}[dr] \ar@{-}[d] & x^4 \ar@{>}[r] \ar@{>}[dr] \ar@{-}[d] & \cdots  &   \\
&  ux \ar@{>}[r]  \ar@{>}[rd] \ar@{-}[d] & ux^2 \ar@{>}[r] \ar@{>}[rd] \ar@{-}[d] & ux^3 \ar@{>}[r] \ar@{>}[dr] \ar@{-}[d] & ux^4 \ar@{>}[r] \ar@{>}[dr] \ar@{-}[d] & \cdots  &   \\
&  u^2x \ar@{>}[r]  \ar@{>}[rd] \ar@{-}[d] &  u^2x^2 \ar@{>}[r] \ar@{>}[rd] \ar@{-}[d] &  \underline{u^2x^3} \ar@{>}[r] \ar@{>}[dr] \ar@{-}[d] &  \underline{u^2 x^4} \ar@{>}[r] \ar@{>}[dr] \ar@{-}[d] & \cdots  &   \\
&  u^3x \ar@{>}[r]  \ar@{>}[rd] \ar@{-}[d] &  u^3x^2 \ar@{>}[r] \ar@{>}[rd] \ar@{-}[d] &  u^3x^3 \ar@{>}[r] \ar@{>}[dr] \ar@{-}[d] &  \underline{u^3 x^4} \ar@{>}[r] \ar@{>}[dr] \ar@{-}[d] & \cdots  &   \\
 &  u^4x \ar@{>}[r]  \ar@{>}[rd] \ar@{-}[d] &  u^4x^2 \ar@{>}[r] \ar@{>}[rd] \ar@{-}[d] &  u^4x^3 \ar@{>}[r] \ar@{>}[dr] \ar@{-}[d] &  \underline{u^4 x^4} \ar@{>}[r] \ar@{>}[dr] \ar@{-}[d] & \cdots  &   \\
 & u^5x \ar@{>}[r] &u^5x^2 \ar@{>}[r] &u^5x^3 \ar@{>}[r] &\underline{u^5x^4}\ar@{>}[r] &\cdots  \ . &
}
$$
\underline{Legend:}
\begin{itemize}
  \item[--] Horizontal arrows: multiplication by $x$.
  \item[--] Oblique  arrows: multiplication by $ux$. Multiplications by other tangibles are not indicated.
  \item[--] Vertical edges: same fiber over  $M$.
  \item[--] Underlined elements: elements in $\mfa_{2,3}$.

\end{itemize}
\end{diagram}

Note for example that
\begin{enumerate}\ealph
  \item $x \cdot u x^3 = ux^4 \notin\mfa_{2,3}$, $ux \cdot u x^3 = u^2 x^4 \in\mfa_{2,3}$,
  \item $x \cdot  x^2 = x^3 \notin\mfa_{2,3}$, $u^2x \cdot x^2 = u^2 x^3 \in\mfa_{2,3}$,
  \item $x \cdot  x^4 = x^5 \notin\mfa_{2,3}$, $u^3x \cdot x^4 = u^3 x^5 \in\mfa_{2,3}$.
\end{enumerate}
This proves that the pairs $(x,ux)$, $(x,u^2 x)$, $(x,u^3 x)$ have different fates, but as is clear from the diagram, $u^3 x$, $u^4 x$, and $u^5 x$  have the same fate.

 Finally we obtain  the supertropical monoids $U'_{i,j}/F$ with  fate distinction. Note that $U'_{i,j}/F$ and $U'_{k,\ell}/F$ are not isomorphic for $(i,j) \neq (k,\ell)$.
For example, for the tangible elements of $U'_{2,3}/F$ we have the following diagram:
\begin{diagram}
 $$
\xymatrix@R=1.5em@C=1.8em{  1  \ar@{>}[r] \ar@{>}[dr]&  x \ar@{>}[r]  \ar@{>}[rd] \ar@{-}[d] & x^2 \ar@{>}[r] \ar@{>}[rd] \ar@{-}[d] & x^3 \ar@{>}[r] \ar@{>}[dr] \ar@{-}[d] & x^4 \ar@{>}[r] \ar@{>}[dr] \ar@{-}[d] & \cdots  &   \\
&  ux \ar@{>}[r] \ar@{-}[d] \ar@{>}[dr]  & ux^2 \ar@{>}[r] \ar@{-}[d] & ux^3 \ar@{>}[r] \ar@{-}[dd]& ux^4 \ar@{>}[r]  & \cdots  &   \\
&  u^2x  \ar@{>}[rd] \ar@{-}[d] \ar@{>}[r] &  u^2x^2  \ar@{>}[rd] \ar@{-}[d] &&&&& &   \\
&  [u^3x]_F \ar@{>}[r]  &  [u^3x^2]_F \ar@{>}[r] & [u^3x^3]_F \ . && \\ }
$$
\end{diagram}

A good hold on the fates of all tangible elements remains reachable, if we replace the equalizer
$\Eq(u^i x^j, c^j)$ by a  fiberwise equalizer obtained as follows. Instead of the monomial~ $u^i x^j$, pick a set $T$ of monomials in $\mfA'$, disjoint from $\{ 1\} \cup \{ u^i x \pipe i \geq 0 \}$, and introduce the the fiberwise equalizer
$$E_T := \Feq(T \cup eT) = \bigvee_{z \in T} \Eq(z,ez),$$
which identifies each monomial $u^i x^j \in T$ with $eu^i x^j= c^j$.
This gives the supertropical monoid
\begin{equation}\label{eq:8.5}
 U'_T : = \mfA'/E_T
\end{equation}
over $M = e\mfA'$.

Let $\mfa_T$ denote the monomial ideal generated by $T$ in $\tT(\mfA')$. The classes $[u^i x^j]_E$ with $u^i x^j \in \mfa_T$ are those classes which contain elements of $M$, namely $c^j$, and so we write $[u^i x^j]_E = c^j$. Every other class $[u^i x^j]_E$  contains only the monomial $u^i x^j$, and so we denote it by $u^ix^j$ itself. The latter are the tangible elements of $U'_T$.

\begin{lem}\label{lem:8.9} Given $u^i x^j \in T$, $\mfa_T$ contains the set
$\{ u^p x^q \pipe p \geq i, q \geq j+1\}$  and is disjoint from
$\{ u^p x^j \pipe p > i\}$.
\end{lem}
\begin{proof}
If $p > i$, $q > j$, then $u^p x^q = u^i x^j \cdot u^{p-i} x^{q-j} \in \mfa_T $, since $u^{p-i} x^{q-j} \in \mfA'$. Also $u^{i} x^{j} \in \mfa_T$. But $u^{p} x^{j} = u^{i} x^{j} u^{p-i} \notin \mfa_T$, since $u^{p-i}  \notin \mfA'$.
\end{proof}

\begin{defn}$ $
\begin{enumerate}\ealph
  \item A monomial $u^i x^j \in \mfA'$ is called a \textbf{border tangible}, if $u^i x^j \notin \mfa_T$ but $u^i x^{j+1} \in \mfa_T$. (N.B. $\mfA' \sm \mfa_T = \tT(U'_T)$.)
  \item    A monomial $u^i x^q \in T$ is called \textbf{basic}, if $u^i x^{q-1}$ is border tangible. Otherwise $u^i x^q$ is said to be \textbf{superfluous} (in $\mfa_T$).

\end{enumerate}
\end{defn}

\begin{prop}\label{prop:8.11} The set of basic elements is the unique minimal set of generators of ~$\mfa_T$.
\end{prop}
\begin{proof}
  This follows from the fact that $u^p x^q$ is superfluous iff there exists $u^i x^j \in T$ with $i \leq p$, $j +1 \leq q$. Thus we can omit all  superfluous elements in the set $T$ of generators of ~$\mfa_T$.
\end{proof}

The following can now be verified in a straightforward way.

\begin{schol}\label{schl:8.12} Let $\{ u^{i_1} x^{j_1}, u^{i_2} x^{j_2}, \dots\} $ be the set of basic elements in $T$, ordered such that $i_1< i_2 <  i_3 < \cdots $. (N.B. obviously all indices $i_k$ are different, since all powers of $x$ are in $\mfA'$.)
  \begin{enumerate} \ealph
    \item The horizontal sequences $\{ u^{i} x, u^{i} x^{2}, \dots\} $ with $i < i_1$ are in $\tT(U'_T)$.  The  horizontal sequences with $i \geq  i_1$  contain  elements of $\mfa_T$, whence each of these  contains a unique border tangible -- the last tangible element of $U'_T$ in  each row.
    \item Every monomial $u^{i_k}x^{j_k-1}$ is border tangible.

    \item If there is a last basic element $u^{i_r}x^{j_ r}$, then for every $i > i_r$ the monomial $u^{i}x^{j_r}$ is border tangible.

        \item Assume again that there is  a last basic element $u^{i_r}x^{j_r}$.  Then the $U'_T$-tangible elements in each column of the row starting with $x, ux, \dots, u^{i_r}x, u^{i_r+1}x $ have different fates, while the $U'_T$-tangibles in each row starting with $u^p x$, $p \geq i_r+1$ all have the same fate.

  \end{enumerate}

\end{schol}
As before we do not give a formal proof, but present an example, inviting the reader to verify this by hand.

\begin{examp}
  Let
  $ T = \{ u^2 x^4,  u^3 x^3,  u^4 x^2,  u^7 x^2\} $, for which all monomials are basic.
  \\ List of tangibles in $U_T$:
  \begin{itemize}
    \item[--] $1, x, x^2, x^3,  x^4, \dots $

    \item[--] $ ux, ux^2, ux^3,  ux^4,  \dots $
    \item[--] $u^2 x, u^2 x^2, u^2x^3$
    \item[--] $u^3 x, u^3 x^2$
    \item[--] $u^4 x$
    \item[--] $u^5 x, u^5 x^2$
    \item[--] $u^6 x, u^6 x^2$
    \item[--] $u^7 x$
    \item[--] $u^8 x, u^8 x^2$
    \item[--] $u^9 x, u^9 x^2$
    \item[--] $\ \vdots \qquad   \vdots$ .

  \end{itemize}
  The monomials in the same column of the row starting with $x, ux, \dots, u^8x$ have different fates, while those in the row starting with  $u^8x, u^9x, \dots$ have the same fate.
\end{examp}

Finally, we look for tyrants in the examples of supertropical monoids just studied, in particular for the \textbf{lonely tyrants}, by which we mean those tangibles which do not have a proper son.

\begin{schol}[List of tyrants]\label{schol:8.14}
$ $
\begin{enumerate}\ealph
  \item In Example \ref{examp:8.5} the only tyrants are $x^{d-1} $ and $u^{d-1}$, and these are lonely.
  \item In the  supertropical monoid $U'_T$ (cf. \eqref{eq:8.5} and Scholium \ref{schl:8.12}) the elements $u^{i_1-1} x^{j_1},$ $ u^{i_1-1} x^{j_1+1}, \dots$ in the row starting with
      $u^{i_1-1} x$ are tyrants, but no one is lonely. The border tangibles in the row starting with
      $u^{p} x$, $p \geq i_r+1$ (cf. Scholium \ref{schl:8.12}.(d)) are lonely tyrants. These are all tyrants in $U'_T$.
  \item In $U_T'/F$ the elements $u^{i_1-1} x^{j_1}, \dots$ are still tyrants, not lonely, but the $F$-class $[u^{i_r-1} x]_F$ and its sons $[u^{i_r-1} x^2]_F, \dots $ are tyrants, the last one is lonely. Remarkably, another new tyrant $u^{i_r-1} x^{j_r -1} = [u^{i_r-1} x^{j_r -1}]_F$ arises, and is not lonely.

      \item Note also that in $U'_T$ and $U'_T/F$ all tangibles are isolated.
\end{enumerate}

\end{schol}

\section{A class of supertropical monoids  with isolated tangibles, which are governed by a lonely tyrant} \label{sec:9}

We introduce another family of supertropical monoids which have fate distinction, and have no zero divisors. In what follows
$$ M = \{ e, c, c^2 , \dots \} \ \dot \cup \ \{ 0 \}  $$
 is the unique totally ordered monoid with $c > e$, and
$\{ e, c, c^2 , \dots \}$ is the cyclic ordered semigroups (all $c^i$ different) together with a specific  absorbing element $0$ as considered in ~\S\ref{sec:8} (cf. \eqref{eq:8.b}).

We start with the ``universal'' supertropical monoid  $\mfA = \mfA(I)$ over $M$  on an indexed set of letters $(t_i \pipe i \in I)$ for any set $I$, where $\tT(\mfA)$ is the free abelian semigroup on $(t_i \pipe i \in I)$ and $e t_i  =c $ for all $i \in I$. Thus the elements of $\tT(\mfA)$ are the monomials in variables $t_i$. We write them as the formal products
\begin{equation}\label{eq:9.1}
  z = \prod_{i \in I } t_i^{\al_i}
\end{equation}
  with exponents $\al_i \in \N_0$ and finitely many $\al_i > 0$. (Read $z = 1$, if all $\al_i =0 $.) Note that~ $ez$ stores the degree of the monomial $z$,
  \begin{equation*}\label{eq:9.2}
    ez = c^{\sum_i \al_i}.
  \end{equation*}
  We choose a monomial
\begin{equation*}\label{eq:9.3}
  \om = \prod_{i \in I } t_i^{r_i} \neq 1
\end{equation*}
and introduce on $\mfA$ the fiberwise equalizer
\begin{equation*}\label{eq:9.4}
  E_\om = \bigvee_{i \in I } \Eq(t_i^{r_i+1}, c^{r_i + 1}).
\end{equation*}

We then obtain the supertropical monoid
\begin{equation*}\label{eq:9.4}
  A_\om = A(I)_\om  := \mfA(I)/E_\om
\end{equation*}
over $M$. This supertropical monoid $A_\om$ has finitely many tangibles, uniquely represented by the monomials $z|\om$, i.e., with  $zw = \om$ for some monomial $w \in \tT(\mfA)$. We denote the $E_\om$-class of every such monomial $z$  by $z$ itself, and the $E_\om$-class of a monomial $z \nmid \om$ by its ghost representative $ez \in M$.

Note that $A_\om$ has no zero divisors, since $M$ has no zero divisors and $E_\om$ is a fiber contraction. It is also evident that every tangible of $A_\om$ is isolated, since every $u \in \tT(A_\om) \sm \{ 1 \}$ has image $eu  > e$ in $M$.

\begin{thm}\label{thm:9.1} $ $
\begin{enumerate}\ealph
  \item For each tangible $z$ of $A_\om$ there exists a unique tangible $z^*$ of $A_\om$ with $z z^* =\om$.
  \item If $z $ and $w$ are tangible of $A_\om$ with $ez = ew $, $z \neq w$, then $z w^* \in M$, $w z^* \in M$, and $z,w $, have different fates. Thus $A_\om$ has fate distinction.

  \item If $z$ is a tangible in $A_\om$, then the sons of $z$ in $A_\om$ are the monomials
  $zw$ with $w | z^*$.
  \item $\om$ is the only lonely tyrant in $A_\om$. The other tyrants of $A_\om$ are the elements $x_i^*$ with $i \in I$, $x_i | \om$.
\end{enumerate}
\end{thm}

\begin{proof} (a): Let $z \in \tT(A_\om)$. Then $z |w $ in $\tT(\mfA)$, and so $zz^* = \om$ for a unique monomial $z^*$ in $\tT(\mfA)$. Since $z^* | \om$, also $z^* \in \tT(A_\om)$.
\pSkip
(b): Let $z ,w \in \tT(A_\om)$, $ez = ew$,  $z \neq w$. Then
$ewz^*= ez z^*$, but $wz^* \neq zz^*= \om$ in $\mfA$. This implies that $wz^*$ is ghost in $A_\om$, since $w$ is the only monomial in $\tT(A_\om)_{e\om}$. For the same reason $zw^*$ is ghost in $A_\om$. From $zz^*\in \tT(A_\om)$, $wz^*\in \tG(A_\om)$, $zw^*\in \tG(A_\om)$, and $ww^*\in \tT(A_\om)$ we infer that $z$ and $w$ have different fates in $A_\om$.
\pSkip
(c): Obvious from $z z^* = \om$.
\pSkip
(d): The claims regarding  tyrants are an immediate consequence of (c).
  \end{proof}

\begin{defn}\label{def:9.2} Given $z \in \tT(A_\om)$, we call $z^*$ the \textbf{complement} of $z$ in $\tT(A_\om), $ or say, that $z^*$ is \textbf{complementary tangible to} $z $ in $A_\om$.
\end{defn}

We regard $A_\om $ as a subset of $\mfA$, using the above notation for the $E_\om$-classes in $\mfA$. Nevertheless, $A_\om$ is not a submonoid of $\mfA$, since the multiplication  in $A_\om$ is different from the multiplication in $\mfA$. Denoting the product of two monomials $z_1, z_2 \in \tT(A_\om)$ by $z_1 \cdot_\om  z_2$, we have
\begin{equation*}\label{eq:9.6}
  z_1 \cdot_\om  z_2 = \left\{ \begin{array}{ll}
                                 z_1 z_2 & \text{if } z_1 z_2 | \om, \\
                                 e z_1 z_2  & \text{else }.
                               \end{array}
  \right.
\end{equation*}

We  look for relations between the monoids  $A_\om \subset \mfA$ for varying $\om \in \tT(\mfA)$.
Given subsets~ $B$ and $C$ of $\mfA$, let $B \cdot C$ denote the set of all product $bc$ with $b \in B$, $c \in C$ in the monoid~ $\mfA$.
\begin{prop}\label{thm:9.3}
  Let $\om_1, \om_2 \in \tT(\mfA)$.
  \begin{enumerate}\ealph
    \item Every $z \in \tT(A_{\om_1 \om_2})$ is a product $z = z_1 z_2$ of monomials $z_1 \in \tT(A_{\om_1})$, $z_2 \in \tT(A_{\om_2})$.
    \item $A_{\om_1 \om_2} = A_{\om_1} \cdot A_{\om_1}.$
  \end{enumerate}
\end{prop}
\begin{proof}
  (a): If $z_1 | \om_1$ and $z_2 | \om_2$, then $z_1 z_2 | \om_1 \om_2$. Conversely, suppose that $z | \om_1 \om_2$, and  write $\om_1 = \prod_{i \in I} t_i^{r_i}$,  $\om_2 = \prod_{i \in I} t_i^{s_i}$, $z = \prod_{i \in I} t_i^{\gm_i}$ (cf. \eqref{eq:9.1}). For a fixed $i \in I$, we have $\gm_i \leq r_i + s_i$. By elementary reasoning there exists a  pair $(\al_i, \bt_i) \in \N_0 \times \N_0$ such that $\al_i \leq r_i$, $\bt_i \leq s_i$, $\al_i + \bt_i =  \gm_i$. Thus $z$ has a  factorization $z = z_1 z_2$ with  $z_1 | \om_1$, $z_2 | \om_2$.
  \pSkip
  (b): $A_{\om_1} \cdot A_{\om_2} = [M \cup \tT(A_{\om_1})] \cdot [M \cup \tT(A_{\om_2})] =
  M \cup [\tT(A_{\om_1}) \cdot  \tT(A_{\om_2})] =
  M \cup \tT(A_{\om_1} \cdot  A_{\om_2}) =
  A_{\om_1\om_2}$.
\end{proof}

\begin{rem}\label{rem:9.4}
Clearly the following hold for
 subsets $A_{\om_1}$ and $A_{\om_2}$ of $\mfA$:
$$ A_{\om_1} \wedge  A_{\om_2} = A_{\om_1} \cap A_{\om_2}$$
and
$$ A_{\om_1} \subset  A_{\om_2} \dss \Leftrightarrow \om_1 | \om_2. $$
\end{rem}

We  present $A_{\om_1} \cdot A_{\om_2}$   as a fiber product.
Given two monomials $z_1, z_2 \in \mfA$, $$z_1 = \prod_{i \in I} t_i^{\al_1}, \qquad   z_2 = \prod_{i \in I} t_i^{\al_2}, $$ we denote the least common multiple of $z_1$
 and $z_2$ as
 \begin{equation*}\label{eq:9.7}
   z_1 \vee z_2 := \prod_{i \in I} t_i^{\max(\al_1, \al_2)},
 \end{equation*}
 and the greatest common divisor as
  \begin{equation*}\label{eq:9.8}
   z_1 \wedge  z_2 := \prod_{i \in I} t_i^{\min(\al_1, \al_2)}.
 \end{equation*}
Clearly
  \begin{equation}\label{eq:9.9}
   z_1 z_2 = (z_1 \wedge  z_2)(z_1 \vee  z_2) .
 \end{equation}

 Given two monomials $\om, \om'  \in \tT(\mfA) $ where $\om| \om'$, we define a projection
\begin{equation*}
  p_{\om',  \om}: A_{\om'}  \Onto  A_\om
\end{equation*} by
\begin{equation}\label{eq:9.12}
  p_{\om',  \om}(z) = \left\{ \begin{array}{ll}
                              z & \text{if } z \in M, \\
                              z \wedge \om  & \text{if } z \in \tT(A_{\om'}).
                            \end{array}
  \right.
\end{equation}

\begin{thm}\label{thm:9.4} $ $
\begin{enumerate}\ealph
  \item Given $\om_1, \om_2 \in \tT(\mfA)$, the monoid $A_{\om_1 \vee \om_2}$ is the fiber product of $A_{\om_1}$ and $A_{\om_2}$ via the projections ($i=1,2$)
$$p_i := p_{\om_1\vee \om_2, \, \om_i}: A_{\om_1\vee \om_2} \To  A_{  \om_i},$$
$$q_i := q_{\om_i, \, \om_1 \wedge  \om_2}: A_{\om_i} \To A_{\om_1 \wedge \om_2}.$$
In short
\begin{equation*}
  A_{\om_1\vee \om_2} = A_{\om_1} \times_{A_{\om_1 \wedge \om_2}} A_{\om_2}.
\end{equation*}

  \item In particular, if  $z_1 \in \tT(A_{\om_1})$, $z_2 \in \tT(A_{\om_2})$, and $z_1 \wedge \om_2 = \om_1 \wedge z_2$, then $z_1 \vee z_2 $ is the unique monomial  $
z  \in \tT(A_{\om_1 \vee \om_2})$ with $z \wedge \om_1 = z_1$, $z \wedge \om_2 = z_2$.
\end{enumerate}
\end{thm}
\begin{proof}
We proceed in  three steps.
\pSkip 1) If $z \in A_{\om_1 \vee \om_2}$, then of course $q_1(p_1(z)) = q_2(p_2(z))$.
\pSkip
2) Let $z_1 \in A_{\om_1}$, $z_2 \in A_{\om_2}$ be given with $q_1(z_1) = q_2(z_2)$. If $z_1 \in M $, then $z_1= q_1(z_1) = q_2(z_2)$, which by \eqref{eq:9.12} implies that $z_1 = z_2$. For the same reason $z_2 \in M $ implies $z_1 = z_2$.
\pSkip
3) There remains the case that $z_1 \in \tT(A_{\om_1})$, $z_2 \in \tT(A_{\om_2})$, i.e.,
$z_1,z_2 \in \tT(\mfA)$, $z_1 | \om_1$, $z_2 | \om_2$, and that $z_1 \wedge \om_2 = \om_1 \wedge z_2$. Now claim (b) holds, since $\tT(\mfA)$ is a distributive lattice with respect to divisibility. For convenience  we give the argument in details. If
$z | (\om_1 \vee \om_2)$, then $z = z \wedge (\om_1 \vee \om_2)  =  (z \wedge \om_1) \vee (z \wedge \om_2)$. If $ z_1 | \om_1$, $z_2 | \om_2$,  and $z_1 \wedge \om_2 = \om_1 \wedge z_2$,
then $(z_1 \vee z_2) \wedge \om_1 = (z_1 \wedge \om_1) \vee (z_2 \wedge \om_1) =
z_1 \vee (z_2 \wedge \om_1) =z_1 \vee (z_1 \wedge \om_2)
= (z_1 \vee z_1) \wedge (z_1 \vee \om_2) = z_1 $. For the same reason $(z_1 \vee z_2) \wedge \om_2 = z_2$.
\end{proof}
\begin{thm}\label{thm:9.5}
  Let $\om_1, \om_2 \in \tT(\mfA)$.
  \begin{enumerate}\ealph
  \item  Then $A_{\om_1} \cdot A_{ \om_2} = A_{\om_1 \wedge \om_2} \cdot A_{\om_1 \vee  \om_2} =  A_{\om_1 \om_2}.$

    \item In particular, a monomial $z \in \tT (\mfA)$ divides $\om_1 \om_2$ iff there exist monomials $u,v$ such that $u | \om_1 \wedge \om_2$, $v | \om_1 \vee \om_2$, and $z = uv.$ Then $u$ and $v$ are uniquely determined by $z$.

    \end{enumerate}
\end{thm}
\begin{proof} (a): Use Theorem \ref{thm:9.4}.(b) and the identity $\om_1\om_2 = (\om_1 \wedge \om_2) \cdot (\om_1\vee \om_2) $, cf. \eqref{eq:9.9}.
\pSkip
(b):  Compare the tangible parts  in (a) and recall Theorem \ref{thm:9.4}.(a).
\end{proof}

\section{Factorization of tangible elements into irreducible elements} \label{sec:10}

Given a supertropical monoid $U$, we say that  an element $x$ of $U$ is \textbf{irreducible}, if $x \in \tT(U)$, but $x \neq 1$, and $x$ cannot be written as a product $uv$ of two elements $u,v \in \tT(U) \sm \{ 1 \} $. Note that in our setting  irreducibility is observed only for tangible elements, and not for all elements $x \neq 1$ as in the standard monoid theory.
The set of all irreducible elements of $U$ is denoted by $\Irr(U)$. We are interested in the case that every $z \in \tT(U)$, $z \neq 1$, is a product of irreducible elements. Then, choosing an indexing,
\begin{equation*}\label{eq:10.1}
  \Irr(U) = \{ x_i \pipe  i \in I \},
\end{equation*}
and we  write every $z \in \tT(U)$ as
\begin{equation}\label{eq:10.2}
  z  = \prod_{i \in I } x_i^{\al_i}
\end{equation}
with $\al_i \in \N_0$, where almost all $\al_i = 0$. ($z=1$ if all $\al_i=0.$) We say that $U$ has \textbf{unique tangible factorization} (abbreviated \textbf{UF}), if each $z \in \tT(U)$ has only one presentation~ \eqref{eq:10.2}, i.e., the exponents $\al_i$ are uniquely determined by $z$.  We want to understand, when~ $U$ has UF, and want to measure in some way, how far $U$ is away from being UF otherwise.

It seems mandatory to assume that the totally ordered set $\tG(U)= M \sm \{0 \} $ is a discrete  semigroup under multiplication, so that it makes sense to speak on irreducible elements in the monoid $\tG(U)$. We will be content with the case that $\tG(U)$ is generated by a single element $c > e$. Thus we assume, that \emph{all considered supertropical monoids  have the ghost ideal}
\begin{equation*}\label{eq:10.3}
  M  = \{ e, c, c^2, \dots  \} \ds{\dot\cup}  \{ 0 \},
\end{equation*}
with $c > e$, and all $c^i$ are different.

In the following we denote the set of all $z \in \tT(U)$ with $ez = d$ for a given $d \in \tG(U)$ by $\tT_d(U)$ instead of $\tT(U)_d$, to avoid a clash with the subscripts $(A_\om, \dots )$ used in \S\ref{sec:9},  which will come up below.
\begin{lem}\label{lem:10.1}
  Assume that $U$ is a supertropical monoid over $M$ with $\tT_e(U) = \{1 \}$. Then every $x \in \tT_c(U)$ is irreducible.
\end{lem}

\begin{proof}
  This follows from the fact that for $u,v \in \tT(U) \sm \{ 1 \}$ we have $eu \geq c $, $ev \geq c$, and so $e(uv) \geq c^2$.
\end{proof}

Later we will frequently meet instances where all irreducible elements are in $\tT_c(U)$ and  in addition $\tT(U)$ is finite.  We now take a brief look at this case.
\begin{lem}\label{lem:10.2}
  Assume that $U$ is a supertropical monoid over $M$, and that every element of $ \tT(U)$ is a product of elements of $\tT_c(U)$. (In particular $\tT_e(U) = \{ 1\}$.)
  \begin{enumerate} \ealph
    \item    Then  $ \tT_c(U) = \Irr(U) = \{ x_i \pipe i \in I\} $.

    \item The set $\tT(U)$ is finite iff $I$ is finite, and for every $i \in I $ there is an exponent $d_i \in \N$ such that $x_i^{d_i} \in \tT(U)$, but $x_i^{d_i +1 } \in M$.
        Then for every monomial  $z = \prod_{i \in I} x_i ^{\al_i} \in \tT(U)$  we have  $0 \leq \al_i \leq d_i$ for each $i \in U$, and so
        $$ |\tT(U)| \leq \prod_{i \in I} (d_i+1). $$

  \end{enumerate}
\end{lem}

\begin{proof}
(a): Evident by Lemma \ref{lem:10.1} and the definition of irreducibility. \pSkip
(b): Now obvious.
\end{proof}

Before delving into details about factorization of tangibles into irreducibles, we mention general features   of the monoids studied here.
\begin{lem}\label{lem:10.3}
  Assume that $U$ is a supertropical monoid over $M$.

  \begin{enumerate} \ealph
    \item Then $U$ is a supertropical semiring.
    \item $U$ has no zero divisors.
    \item If $\tT_e(U) = \{1 \}$, then every tangible element of $U$ is isolated.
  \end{enumerate}
\end{lem}

\begin{proof}
(a): The criterion in \cite[Theorem 1.2]{IKR4}, that $U$ is a supertropical semiring, holds for a trivial reason: If $0 < ex < ey $ for $x,y,z \in U$, then $0 < exz < eyz$ for $z \neq 0$, since $M \sm \{ 0 \}$ is cancellative under multiplication, while $yz= eyz = 0 $ when $z= 0$.
\pSkip
(b), (c): Recall the arguments preceding to Theorem \ref{thm:9.1}.
\end{proof}

Essentially in what follows, $U$ is assumed to be a supertropical monoid with $eU = M$, having the following properties:
\begin{enumerate} \dispace
  \item [(A)] $ \Irr(U) = \tT_c(U) =  \{ x_i \pipe i \in I\} $.
   \item [(B)] $\tT_e(U) = \{ 1\} $; every $z \in \tT(U)$ is a product of elements in $\tT_c(U)$.
\end{enumerate}
Monoids with these  properties abound in the category of supertropical monoids over $M$, and arise as follows.

\begin{rem}\label{rem:10.4}
Given a supertropical monoid over $M$, let $\langle \tT_c(U)\rangle_\tT$ denote the set of all finite products of $\tT_c(U)$ which are tangible.
The subset
$$ U' := \{ 1 \} \ds{\dot \cup } \langle \tT_c(U) \rangle_\tT \ds{\dot \cup } M $$
is a submonoid of $U$ with  $eU' = M$, admitting properties (A) and~ (B).
\end{rem}

\emph{\textbf{Starting from  now $U$ is always a supertropical monoid over $M$ that  has properties (A) and ~(B).}}

\pSkip

 To analyze the factorization of tangibles in $U$
 we use the supertropical monoid $\mfA := \mfA(I)$ introduced in \S\ref{sec:9}, where now we identify the letters $t_i \in \tT_c(U)$ with the elements $x_i \in \tT_c(U)$. Note that $e \mfA = e U  = M$.
 We have the unique factorization $$ \pi_U: \mfA \twoheadrightarrow U$$
 over $M$ which maps every $x_i \in   \tT_c(\mfA)$ to $x_i \in \tT_c(U)$ and every $z \in M $ to itself. Let $E_U$ denote the MFCE-relation on $M$ induced by $\pi_U$, i.e., $E_U$ is the set of pairs $(z_1, z_2) \in \mfA \times \mfA$ with $\pi_U(z_1) = \pi_U(z_2)$. We have a natural fiber contraction
 $$ U = \mfA / E_U. $$

 When the $E_U$-class $[z]_U = \pi_U(x)$ of a monomial $z \in \tT(\mfA)$ is tangible, this element of $U$ will often be denoted by $\olz$, while if $[z]_U$ is ghost it will be denoted by its unique representative $ez$ in $M$.  Observe that in this notation a product $\olz_1 \cdot \olz_2$ can be different from $[z_1  z_2] _U$. More precisely,  $\olz_1 \cdot \olz_2 = \overline{ z_1  z_2}$, if $\olz_1 \cdot \olz_2 \in \tT(U)$, and  $\olz_1 \cdot \olz_2 = e z_1 z_2$ otherwise. \pSkip

The following  evident fact is crucial  for our approach to unique factorization.
 \begin{rem}\label{rem:10.5}
An element $\zt$ of $\tT(U)$ has a unique factorization $\zt = \prod_{i \in I} x_i ^{\al_i}$ iff there is only one $ z \in \tT(\mfA)$ with $\zt = \olz$. (Recall that we regard every $x_i$ as element of $\tT_c(\mfA)$ and of $\tT_c(U)$.)
\end{rem}

\begin{defn}\label{def:10.6}$ $
\begin{enumerate}
  \ealph
  \item We call a fiber contraction $W$ of $\mfA$ over $M$ with $E_W \subset E_U$ a \textbf{cover of} $U$.
  \item We then have the projection $\pi_{W,U} : W \onto U$ over $M$ with $\pi_{W,U} \circ \pi_W = \pi_U$.
  \item Given covers $W_1$ and $W_2$ of $U$, $W_2$ is called a \textbf{subcover} of $W_1$, if
  $E_{W_1} \subset E_{W_2} \subset E_{U}$. Then
  $\pi_{W_1,U} = \pi_{W_2,U} \circ \pi_{W_1,W_2}$.
\end{enumerate}
\end{defn}
\noindent N.B. We do not need to say ``cover of $U$ in $\mfA$'', since the monoid $\mfA$ is uniquely determined by $U$ (up to unique isomorphism over $U$) due to the  properties (A) and (B).
\pSkip

We denote the set of all covers of $U$ by $\Cov(U)$. It is partially ordered by the relation ``subcovers''.
\begin{schol}\label{schl:10.7}
  The set of all MFCE-relations on $\mfA$, ordered by inclusion, is a complete lattice $\MFC(\mfA)$, as we know for long \cite[\S6]{IKR1}, \cite[\S1]{IKR4}. Its subset $\MFC(\mfA, U)$ consisting of the relations containing $E_U$ is again a complete lattice. Consequently, the set of covers of $U$ is a complete lattice anti-isomorphic to $\MFC(\mfA,U)$. For $W_1,W_2 \in \Cov(U)$  we have
  $W_1 \leq W_2 \Iff E_{W_1} \supset E_{W_2}$, $E_{W_1 \vee W_2} = E_{W_1} \cap E_{W_2} = E_{W_1} \wedge E_{W_2}$, but $E_{W_1 \wedge W_2} = E_{W_1} \vee E_{W_2}$ may be larger than the set $E_{W_1} \cup E_{W_2}$. Note also that all MFCE-relations on~ $\mfA$ are subsets of $\mfA\times_M \mfA \subset \mfA \times \mfA$.
\end{schol}

\begin{notation} Every equivalence relation $E$ on $\mfA$ is symmetric with respect to the switching  map $s: \mfA \times \mfA \to \mfA \times \mfA$, $s(a,b) = (b,a)$. To ease notation we abbreviate for any  $Z \subset \mfA \times \mfA$  the set $Z \cup s(Z)$ by $Z^\sig$.

\end{notation}

We introduce on $\mfA$ two MFCE-relations, which are coarser than $E_U$, as follows
\begin{equation}\label{eq:10.4}
\begin{array}{ll}
  \tlE_U  := &  \{ (z_1, z_2) \pipe z_1, z_2 \in \tT(\mfA) , \pi_U(z_1) = \pi_U(z_2) \in M \} \\[1mm]
  & \  \dot \cup \ \{ (z,ez) \pipe z\in \tT(\mfA), \pi_U(z) \in M \}^\sig  \\[1mm]
  & \ \dot \cup \ \diag(M),
\end{array}
  \end{equation}
\begin{equation*}\label{eq:10.6}
\begin{array}{ll}
  E'_U  & :=   \big \{ (z_1, z_2) \pipe z_1, z_2 \in \tT(\mfA) , \pi_U(z_1) 
  = \pi_U(z_2) \in \tT(U) \big \}  \ \dot \cup \ \diag(M).
\end{array}
  \end{equation*}
It is  obvious that
\begin{equation}\label{eq:10.6}
  \tlE_U \cup E'_U = E_U, \qquad \tlE_U \cap E'_U = \diag(\mfA).
  \end{equation}
From these relations we obtain two covers
\begin{equation}\label{eq:10.7}
\begin{array}{ll}
 \tlU := \mfA /  \tlE_U, \qquad U' := \mfA / E'_U
\end{array}
  \end{equation}
of $U$.

\newpage
\begin{thm}\label{thm:10.8} $ $
  \begin{enumerate}\ealph
    \item $\tlU \wedge U' = U$ and $\tlU \vee U' = \mfA$.
    \item  $\tlU$ is the unique minimal cover of $U$ with UF.
    \item $U$ has UF iff $\tlU = U$ iff $U' = \mfA$.
    \item $U'$ is the maximal cover of $U$ such that $\pi_{U',U}(\tT(U'))= \tT(U)$. In other words, $U'$ is the maximal cover of $U$ such that $U$ is a tangible fiber contraction of $U'$ (cf. Definition \ref{defn1.3}.c).

  \end{enumerate}
\end{thm}
\begin{proof}
This is now evident from \eqref{eq:10.6}, Scholium \ref{schl:10.7},  and Remark \ref{rem:10.5}.
\end{proof}

We now  take a look at the binary equalizers $\Eq(z,w)$ on $\mfA$ with $ez = ew$. We are in a favorable situation, since  $\tT(\mfA)$ is closed under multiplication, and the monoid $\tT(\mfA)$ is cancellative. Thus we do not need to resort to  the description of equalizers in terms of  paths in \cite[\S6]{IK}.

\begin{prop}$ $
\begin{enumerate}
  \ealph
  \item Let $z_1,z_2 \in \tT(\mfA)$ with $z_1 \neq z_2$, $ez_1 = e z_2$. Then
  \begin{equation}\label{eq:10.8}
    \Eq(z_1, z_2) = \{ (z_1 u , z_2 u ) \pipe u \in \tT(\mfA) \} \ds{ \dot \cup} \diag(\mfA) .
  \end{equation}
  \item
  Let $z \in \tT(\mfA)$. Then
  \begin{equation}\label{eq:10.9}
    \Eq(z, ez) = \{ (z u , e z u) \pipe u \in \tT(\mfA)\}^\sig \ds{ \dot \cup} \diag(\mfA) .
  \end{equation}
 \end{enumerate}
\end{prop}
\begin{proof}
On the right hand side in both cases  there appears a multiplicative equivalence relation on $\mfA$.
To see this it suffices for \eqref{eq:10.8} to use the fact that $\tT(\mfA)$ is closed under multiplication. For \eqref{eq:10.9} we also need the  fact that $\tT(\mfA)$ is cancellative. It is then evident, that these multiplicative equivalences are the smallest  ones containing
 $(z_1,z_2)$ and  $(z,ez)$ respectively.
  \end{proof}
Given a point $\zt \in \tT(U)$, we study covers $W$ of $U$ such that the elements in the fiber
$\pi^{-1}_{W,U}(\zt)$ have a ``better'' factorization into irreducibles than $\zt$. This can be measured by various tools. We first employ binary equalizers as in~ \eqref{eq:10.8}.

\begin{defn}\label{def:10.11} $ $
\begin{enumerate}\ealph
  \item We call a triple $(z_1, u, z_2)$ with $z_1, z_2, u \in \tT(\mfA)$ and $ez_1 = ez_2$ a \textbf{confluence in} ~$U$, if $\overline{z_1 u } = \overline{z_2 u }$, and then say more precisely, that $(z_1, u, z_2)$  is \textbf{confluent at} $\zt = \overline{z_1 u}$.
  \item In this situation we define a cover  $\Con(z_1, u, z_2)$ of $U$ as follows (cf. \eqref{eq:10.8}).
      $$ \Con(z_1, u, z_2) := \mfA/ \Eq(z_1 u , z_2 u ).$$
  \item We call a confluence $(z_1, u, z_2)$ in $U$ \textbf{initial}, if for any monomial $v\| u$ (i.e., $v|u$ and $v \neq u$) the triple $(z_1, v, z_2)$ is not confluent  in $U$.
\end{enumerate}
\end{defn}

\begin{thm}\label{thm:10.12} $ $
\begin{enumerate}\ealph
  \item $\Con(z_1, u, z_2)$ is the unique maximal cover $W$ of $U$ such that
  $$ \zt = \pi^{-1}_{W,U}(z_1 u )  =\pi^{-1}_{W,U}(z_2 u ). $$
  \item If $w \in \tT(\mfA)$ is a monomial such that $\zt \cdot \brw$ is still tangible in $U$, then $\Con(z_1, uw, z_2)$ is a subcover of $\Con(z_1, u, z_2)$ in $\Cov(U)$.
  \item For any confluence $(z_1, u, z_2)$ in $U$ there is a unique least  initial confluence
  $(z_1, u_0, z_2)$ in $U$ with $u_0| u$. The monoid $\Con(z_1, u_0, z_2)$ is the  maximal cover with $\zt \in W_0$, admitting $\Con(z_1, u, z_2)$ as a subcover.
\end{enumerate}

\end{thm}

\begin{proof}
  (a): A direct consequence of the definition of the equalizer $\Eq(z_1u, z_2 u)$ and the fact that $\Eq(z_1u, z_2 u)  \subset E_U.$ \pSkip
  (b): Obvious, since $\Eq(z_1u, z_2 u) \subset \Eq(z_1u w, z_2 uw)$.
  \pSkip
  (c): Now clear, since $\Con(z_1, u,  z_2 )$ is the \emph{unique maximal} cover of $U$ with $\pi^{-1}_{W,U}(z_1 u )  =\pi^{-1}_{W,U}(z_2 u )$.
\end{proof}

We turn to more general equalizers than the binary ones. As before, $\zt$  is a tangible element of the supertropical monoid $U$ over $M$.

\begin{defn}\label{def:10.13} $ $
\begin{enumerate}\ealph
  \item  We say that a cover $W$ of $U$ \textbf{splits $\zt$ totally}, in short ``$W$ splits $\zt$'', if $\ipi_U(\zt) \subset~ W$.

  \item We say that a cover $W$ of $U$ is \textbf{alien} to $\zt$, if  $\ipi_U(\zt) \cap W =  \emptyset$.
\end{enumerate}
\end{defn}
\noindent (a) means that all monomials  $z$ with $\olz = \zt$ are already contained in the monoid $W$.
In other words, all possibilities to factor $\zt$ into  irreducibles of $U$ occur already  in $\tT(W)$ instead of  $\tT(U)$. (N.B. We tacitly assume that $\zt \neq 1$. For $\zt = 1$ our claims are empty.)
For any cover $W$ of $U$ we have $|\ipi_W(\zt) | \leq |\ipi_U(\zt)|$ with equality iff $W$ splits $\zt$ totally and $\ipi_W(\zt)  = \emptyset$ iff $W$ is alien to $\zt$.
Note that, if $W$ is alien to $\zt$,  then certainly $\zt \notin \mfA$.

\begin{thm}\label{thm:10.14} Assume that $\zt \notin \mfA$.  Let $ S := \ipi_U(\zt) \subset \tT(\mfA).$
\begin{enumerate}\ealph
\item There is a unique minimal cover $W$ of $U$ which splits $\zt$ totally. It is the intersection of all covers $W$ of $U$, which contain the set $S$.

  \item Every cover $W$ of $U$, which is alien to $\zt$, is contained in a maximal such cover. The intersection of all maximal covers alien to $\zt$ is the fiber contraction
  \begin{equation*}\label{eq:10.10}
    U_\zt := \mfA / \Eq(\ipi_U(\zt))
  \end{equation*}
  by the MFCE-relation  $\Eq(S). $

\item  If $W$ is a cover of $U$ alien to $\zt$, then $W \vee U_\zt$ is again alien to $\zt.$

  \item  For finite $S = \{ z_1, \dots, z_n \} $, write $z_i = z_i^\circ u $ with $z^\circ _1 \wedge \cdots  \wedge z^\circ _n = 1$, $u = \gcd(z_1, \dots, z_n)$ in $\tT(\mfA)$. Then
      $$ U_\zt = \bigwedge_{i <j } \Con(z^\circ_i, u, z^\circ_j ).  $$
  \item In particular this holds when $\tT_c(U)$ is finite.
\end{enumerate}
\end{thm}
\begin{proof} (a):  Evident from the definition of total splitting.
\pSkip
(b): The first assertion is plain by Zorn's Lemma. To prove the second assertion, we  choose an indexing
$S = \{ z_i \pipe i \in I\} $  with $I$ totally ordered. Given a cover $W$, let $E_W$ denote the MFCE-relation on $\mfA$ over $M$ such that  $W = \mfA/ E_W$. Then $W$ is alien to $\zt$ iff for every $i \in I$ there is an index $j$ with $(z_i,  z_j) \in E_W$, since precisely then $\pi_W(z_i)$ is not a monomial for each $i$.
This means, that a cover $W$ of $U$ is a maximal cover alien to $\zt$ iff $W$  is an intersection of covers $\mfA/ \Eq(z_i, z_j)$, where every element of $\ipi_U(\zt)$ shows up as $z_i$ or $z_j$. It follows that $U_\zt = \mfA/ \Eq(S)$ is the intersection of all these covers, since
$\Eq(S) = \bigvee_{i < j } \Eq(z_i, z_j).$
\pSkip
(c): If $W$ is a cover of $U$ alien to $\zt$, then there is a maximal such cover $W' \supset W$.
Since~ $U_\zt$ and $W$ are both subsets of $W'$,
also $U_\zt \vee W$ is a subcover of $W'$, and so $U_\zt \vee W$ is alien to $\zt$.
\pSkip
(d): Follows from the  fact ,that
$ \Eq(S) = \bigvee_{i < j } \Eq(z_i, z_j),$ and elementary reasoning in the monoid  $\tT(U)$
\pSkip
(e):   Since $U$ has the properties A and B, we know that finiteness of $\tT_c(U)$ implies finiteness of $\ipi_U(\zt)$ for every $\zt \in \tT(U)$, cf. Lemma \ref{lem:10.2}.
\end{proof}

More generally, we look for covers of $U$ which split $\zt$ ``partially''.
\begin{defn}\label{def:10.15}
Let $\zt \in \tT_d(U)$, $d>e$.
We say that a cover $W$ of $U$ is a  \textbf{minimal partial splitting cover for} $\zt$, if for any
$W' \in \Cov(U)$
\begin{equation}\label{eq:10.11}
  W \cap \ipi_U(\zt) = W' \cap \ipi_U(\zt) \dss{\Rightarrow} W \subset W'.
\end{equation}
\end{defn}
Note that then more generally the following holds: For all $W' \in \Cov(U)$
\begin{equation*}\label{eq:10.12}
  W \cap \ipi_U(\zt) \subset  W' \cap \ipi_U(\zt) \dss{\Rightarrow} W \subset W'.
\end{equation*}
Indeed, $W \cap \ipi_U(\zt)$  is contained in $W' \cap \ipi_U(\zt)$  iff
\begin{equation*}
  W \cap \ipi_U(\zt)  =  (W \cap W') \cap \ipi_U(\zt),
\end{equation*}
and then \eqref{eq:10.11} implies that $W \subset W \cap W'$, which means that $W \subset W'$.

Partial splitting of $\zt$ is governed by a family $P\Sig_U(\zt)$  of subsets of $\ipi_U(\zt) \subset \tT(\mfA)$,  defined as follows:
\begin{equation*}\label{eq:10.13}P\Sig_U(\zt) = \{ S \subset \ipi_U(\zt) \pipe \exists W \in \Cov(U) :
 \ipi_U(\zt) \cap W = S\}.
\end{equation*}

\begin{thm}\label{thm:10.16} $ $
\begin{enumerate}\ealph
  \item For every $S \in P\Sig_U(\zt)$ there is a unique minimal cover $W_S$ of $U$ such that
      \begin{equation*}\label{eq:10.14}
        S = \ipi_U(\zt) \cap W_S.
      \end{equation*}
  \item
  $W_S$ is the unique minimal cover of $U$ that contains $S$.
\end{enumerate}
\end{thm}

\begin{proof}
(a): The family of all covers $W$ of $U$ with $\ipi_U(\zt) \cap W = S$ is closed under intersections. Thus the intersection $W_S$ of all these covers is the minimal one cutting out $S$ in~ $\ipi_U(\zt)$.
\pSkip
(b): Clearly, $S \subset W_S$. Let $S \subset W'$. Then $\ipi_U(\zt) \cap W_S \cap W' = S \cap W' = S$. Due to the minimality of $W_S$, we have $W_S \subset W_S \cap W' $, and so $W_S \subset W'$.
\end{proof}

\begin{rems}\label{rem:10.17} The following is now obvious.
\begin{enumerate}\ealph
  \item Let $S_1, S_2 \in P\Sig_U(\zt)$. Then
  $$ W_{S_1} \subset W_{S_2} \dss{\Rightarrow} S_1 \subset S_2. $$
  \item  Assume that $\{ S_i\pipe i \in  I\}$ is a family of subsets of $P\Sig_U(\zt)$. Then
  $ S := \bigcap_{i \in I } S_i \in P\Sig_U(\zt),$ and
  $W_S = \bigcap_{i \in I} W_{S_i}.$
  \item
  In particular, $ P\Sig_U(\zt)$ contains a  smallest set $S_0$, namely the intersection of all $S \in  P\Sig_U(\zt)$.
  %
  \item $S_0 = \emptyset$ iff there exist covers of $U$ alien to $\zt$. Then $W_{\emptyset} = U.$
\end{enumerate}
\end{rems}

\begin{thm}\label{thm:10.18} The  minimal partial splitting covers $W_S$ for  any $\zt \in \tT(U)$ are subcovers of the cover $\tlU$ of $U$ defined above, cf. \eqref{eq:10.7}.
\end{thm}

\begin{proof}
This is a consequence of Theorem \ref{thm:10.8}.b, which states that every monomial $z \in U \cap \tT(\mfA)$ is contained in $\tlU$.
\end{proof}
The partial  splitting concept can be extended without further effort as follows.
Given a family $Z = \{ \zt_j \pipe j \in J\} $ of tangible elements of $U$, i.e., a subset of $\tT(U)$, we look for a cover $W$ of $U$ which splits each $\zt_j$ partially, and then say that $W$ \textbf{splits $Z$ partially}. In the same vein, if $W$ splits each $\zt_j$ totally, we say that $W$ \textbf{splits $Z$ totally}, and if $W$ is alien to each $\zt_j$, we say that $W$ is \textbf{alien} to $Z$.

We proceed as follows. Let
\begin{equation*}\label{eq:10.15}P\Sig_U(Z) = \{ S \subset \ipi_U(Z) \pipe \exists W \in \Cov(U) \text{ such that }
 \ipi_U(Z) \cap W = S\}.
\end{equation*}

\begin{thm}\label{thm:10.19}
$ $
\begin{enumerate} \ealph
  \item Given $S \in P \Sig_U(Z)$, the set of covers $W$ of $U$ with $\ipi_U(Z) \cap W = S$ is closed under intersections. Thus there exists a unique minimal cover $W_S$ of $U$ with
      \begin{equation*}\label{eq:10.16}
        S = \ipi_U(Z) \cap W_S.
      \end{equation*}
The cover $W_S$ is called a \textbf{{minimal partial splitting cover of $U$ for $Z$}}, more precisely the minimal splitting cover for $Z$ and $S$.
\item $W_S$ is the unique minimal  cover of $U$ that contains $S$.  It is a subcover of $\tlU$.

\item  Every cover $W$ of $U$, which is alien to $Z$, is contained in a maximal such cover.
The intersection of all maximal covers alien to $Z$ is the tangible  fiber contraction
  \begin{equation*}\label{eq:10.10}
    U_Z := \mfA / \Eq(\ipi_U(Z))
  \end{equation*}
  of $\mfA$ by the equalizer $\Feq(\ipi_U(Z))$.

  \item  If $W$ is a cover of $U$ alien to $Z$, then $W \vee U_Z$ is again alien to $Z.$
\item Remarks \ref{rem:10.17} remain valid for  $P \Sig_U(Z)$ instead of $P \Sig_U(\zt)$.
\end{enumerate}
\end{thm}
\begin{proof}
These facts can be obtained in a straightforward way by using arguments identical or analogous to those used in the proof of Theorems \ref{thm:10.14}, \ref{thm:10.16}, and \ref{thm:10.18}.
\end{proof}

Note that Theorem \ref{thm:10.19} is \emph{not a trivial consequence} of Theorems \ref{thm:10.14}, \ref{thm:10.16}, and~ \ref{thm:10.18}, since it is by no mean clear, how to obtain a minimal partial splitting cover $W$ for $S$ from the minimal partial splitting covers $W_j$  for the sets
$S_j := S \cap \ipi_U(\zt_j)$, although $S$  is the disjoint union of the sets $S_j$.

\section{Partitioned splitting covers}\label{sec:11}

We use the setting developed in \S\ref{sec:10}, in which  $U$ is a supertropical monoid with ghost ideal
$M = \{ e, c , c^2, \dots \} $, where $c > e$, having properties A and B.
Given a tangible element~ $\zt$ of $U$, we study properties of the set $\ipi_U(\zt)$ relevant for factorization of $\zt$ into irreducibles and their interplay with certain covers of $U$.

As common, a partition $\Pi = \{ X_j \pipe \jJ \} $  of a set $X$ means a decomposition of $X$ into disjoint subsets $X_j$. In other words, $\Pi$ is the set of equivalence relations on the set $X$. The set  of all partitions of $X$ carries a natural partial ordering. We say that \textbf{a partition $\Pi'$ is a refinement of a partition $\Pi$ of $X$}, and write $\Pi \leq \Pi'$, if every set in $\Pi $ is a union of sets in $\Pi'$. This means that the equivalence relation given by $\Pi'$ is finer than the
equivalence relation given by $\Pi$.
\begin{defn}\label{def:11.1} $ $
\begin{enumerate}\ealph
  \item Let $\Pi = \{ S_j \pipe j \in J \} $ be a partition of the set $\ipi_U(\zt)$. We say that a cover $W$ of~ $U$ is \textbf{adapted to} $\Pi$, if the projections $\pi_W(S_j)$ are one-point sets $\{ \zt_j\} \subset \tT(W)$ with $\zt_j \neq \zt_k$ for $j \neq k$. In other terms, there exists a subset $\{ \zt_j \pipe   j \in  J \} $ of $\tT(W)$, such that $\zt_j \neq \zt_k$ for $j \neq k$ and
      \begin{equation}\label{eq:11.1}
        \Pi = \{ \ipi_W(\zt_j) \pipe j \in J \}.
      \end{equation}
            We denote the set of these covers $W$ by $\Sig_U(\zt)$.

  \item If a cover $W$ of $U$ is adapted to $\Pi$, and so \eqref{eq:11.1} holds with different tangibles $\zt_i $ of~ $W$, we say that $\Pi$ is \textbf{induced} by $W$.
\end{enumerate}

\end{defn}

\begin{thm}\label{thm:11.2} Given a partition
$\Pi = \{ S_j \pipe j \in J \} \in \Sig_U(\zt)$, there  a \textbf{unique minimal cover} $W_\Pi$ of $U$ adapted to $\Pi$, namely
\begin{equation}\label{eq:11.2}
W_\Pi := \bigvee_{\jJ} \mfA/\Eq(S_j).
\end{equation}
\end{thm}
\begin{proof}
  Picking a cover $W'$ of $U$ adapted to $\Pi$ (which exists by assumption), we have a set
  $\{ \zt'_j \pipe \jJ \} \subset \tT(W')$ with $\ipi_U(\zt'_j) = S_j$ for every $\jJ$, in particular
  $\zt'_j \neq \zt'_k $ for every $j \neq k$. Set
  $$ W_j : = \mfA / \Eq(S_j).$$
  Recall that $\Eq(S_j)$ is the finest  tangible MFCE-relation on $\mfA$ over $M$ that identifies $S_j$ into one point, now named $\zt_{j o} \in \tT(W_j)$.   Since $\pi_{W'}(S_j) = \{ \zt'_j\} $ is a singleton in $\tT(W')$, we have $W_j \subset W'$ for every $\jJ$. Thus
  $$ W := \bigvee_\jJ W_j \subset W'.$$

For every $\jJ$, let $\zt_j := \pi_{W',W}(\zt'_j).$ Then $\pi_{W,W_j}(\zt_j) = \zt_{jo}$. Suppose that there exist two indices $j \neq k$ such that $\zt_j = \zt_k$, then
$$  \zt'_k = \ipi_{W',W}(\zt_j) = \ipi_{W',W_j}(\zt_{jo}),$$
and so
$$  S_k = \ipi_{W'}(\zt'_k) \subset \ipi_{W'} (\ipi_{W',W_j}(\zt_{jo}))
=\ipi_{W_j}(\zt_{jo}) = S_j.$$
But $S_j \cap S_k = \emptyset.$ We conclude that $\zt_j \neq \zt_k$ for $j \neq k$. Thus $W$ is adapted to $\Pi $ and is contained in every such cover $\Pi'$.
\end{proof}
%

We call the covers \eqref{eq:11.2} obtained in this way the \textbf{minimal partitioned  splitting covers} of $U$ for ~ $\zt$. More precisely, we call $W_\Pi$ the \textbf{minimal partitioned splitting cover} of $U$ for $\zt$ \textbf{inducing}~ $\Pi$.
Since a partition $\Pi$ of $\ipi_U(\zt)$ is uniquely determined by any cover $W$ adapted to $\Pi$, cf. ~ \eqref{eq:11.1}, it is evident that we have a \emph{natural bijection $\Pi \mapsto W_\Pi$ from  $\Sig_U(\zt)$ to  the set of minimal partitioned splitting covers for $\zt$}.

\begin{thm}\label{thm:11.3}
Assume that $W_1$ and $W_2$ are covers of $U$ such that $W_1 \subset W_2$, and that $\Pi_1$ and $\Pi_2$ are partitions of $\ipi_U(\zt)$ adapted to $W_1$ and $W_2$. Then $\Pi_2$ is a refinement of~ $\Pi_1$, i.e.,
$\Pi_1 \leq \Pi_2.$
\end{thm}
\begin{proof}
Let $\Pi_1 = \{ S'_j \pipe j \in J' \} $ and  $\Pi_2 = \{ S''_k \pipe k \in J'' \} $. Since $W_1$ and $W_2$ are adapted to $\Pi_1$ and $\Pi_2$, we have  $S'_j =\ipi_{W_1}(\zt'_j)$, $S''_k =\ipi_{W_2}(\zt''_k)$ with
element $\zt'_j \in \tT(W_1)$ and $\zt''_k \in \tT(W_2)$.
Pick an index $j \in J'$,  and then an index $k \in J''$
with $S'_j \cap S''_k \neq \emptyset$. We verify that $ S''_k \subset S'_j$, and then will done.
We have $\zt'_j = \pi_{W_1}(S'_j \cap S''_k)$ and  $\zt''_k = \pi_{W_2}(S'_j \cap S''_k)$. Thus
$\pi_{W_2, W_1}(\zt''_k) = \zt'_j$, from which we conclude that
$$ S''_k = \ipi_{W_2}(\zt''_k)  \subset \ipi_{W_1}(\zt'_j) = S'_j. $$
\vskip -7mm
\end{proof}
We are ready for the central  result of this section.
\begin{thm}\label{thm:11.4}
Given $\zt \in \tT(U)$, assume that $\Pi_1$ and $\Pi_2$ are two partitions of $\ipi_U(\zt)$, and that $W_1$ and $W_2$ are the minimal splitting covers of $U$ for $\zt$, inducing $\Pi_1$ and $\Pi_2$ respectively.
Then
\begin{equation}\label{eq:11.3}
  \Pi_1 \leq \Pi_2 \dss \Leftrightarrow W_1 \subset W_2.
\end{equation}
Thus the bijection $\Pi \mapsto W_\Pi $ from $\Sig_U(\zt)$ to the set of minimal splitting covers of $U$ for $\zt$ is an isomorphism of partially ordered sets.
\end{thm}
\begin{proof}
The implication $(\Leftarrow)$ in \eqref{eq:11.3} is a consequence of Theorem \ref{thm:11.3}. To prove  $(\Rightarrow)$,  assume that $\Pi_1 \leq \Pi_2$.
Let $\Pi_1 = \{ S'_j \pipe j \in J' \} $ and  $\Pi_2 = \{ S''_k \pipe k \in J'' \} $. For  given $j \in J'$ we have
$$ S'_j  =  \bigcup_{k \in J''_j} S''_k$$
for some subset $J''_j$ of $J''$.
Then $\Eq(S''_k) \subset \Eq(S'_j)$ for every $j \in J'_j$, since  $\Eq(S''_k)$ is the finest MFCE-relation over $M$ that identifies $S''_k$ into a point. We infer that
$$ \bigvee_{k \in J''_j} \Eq(S''_k) \subset  \Eq(S'_j), $$
and so, since $\bigcup_{j \in J'} J''_j = J''$,
$$ \bigvee_{k \in J''} \Eq(S''_k) \subset \bigvee_{j \in J'} \Eq(S'_j).$$
By Theorem \ref{thm:11.2}, this means that $W_1 \subset W_2$.
%
\end{proof}

\pSkip
As an application of Theorem \ref{thm:11.4} we describe how any family of minimal partitioned  splitting  covers of~ $U$ for $\zt$
can be assembled into one such splitting cover.
Let $(\Pi_\lm \pipe \lm \in \Lm)$ be a family of partitions of $\ipi_U(\zt)$, which all are splitting patterns  of $\zt$, $\Pi_\lm = (S_{\lm,j} \pipe j \in J_\lm )$.  Assume that $W_\lm := W_{\Pi_\lm}$ is the minimal splitting cover of $U$ for $\zt$ inducing $\Pi_\lm$. Thus
\begin{equation*}\label{eq:11.4}
  S_{\lm,j} = \ipi_{W_\lm}(\zt_{\lm,j}) \qquad (\jJ)
\end{equation*}
with $\zt_{\lm,j } \in \tT(W_\lm)$, $\pi_U(\zt_{\lm,j}) = \zt. $ The minimal common refinement $\Pi_0 = \bigvee_{\lm \in \Lm} \Pi_\lm $ of all partitions  $\Pi_\lm $  clearly can be  described as follows:
\begin{equation*}\label{eq:11.5}
  \Pi_0 = \{ S_{0,k} \pipe k \in J_0\},
\end{equation*}
where  $J_0$ is the set of families
$(j_\lm \pipe \lm \in \Lm )$ with  $\bigcap_{\lm \in \Lm} S_{\lm,j_\lm } \neq \emptyset$,  and
\begin{equation*}\label{eq:11.6}
S_{0,k} = \bigcap_{\lm \in \Lm} S_{\lm,j_\lm }.
\end{equation*}
In short, $\Pi_0$ is the set of all nonempty intersections of members of the partitions $\Pi_\lm$, where each $\Pi_\lm$ contributes one member. By Theorem \ref{thm:11.4} we know that $W_0 \in  \Sig_U(\zt)$, and that
\begin{equation*}\label{eq:11.7}
W_{0} = W_{\Pi_0 } = \bigvee_{\lm \in \Lm} W_{\lm }
\end{equation*}
is the minimal partitioned splitting cover of $U$ for $\zt$ inducing $\Pi_0$.
We have
\begin{equation*}\label{eq:11.8}
S_{0,k} = \ipi_{W_0}(\zt_{0,k})
\end{equation*}
with unique $\zt_{0,k} \in \tT(W_0)$.
It follows that
\begin{equation*}\label{eq:11.9}
\pi_{W_0,W_\lm}(\zt_{0,k}) = \zt_{\lm,j_\lm}.
\end{equation*}
Instead of an (easy) formal proof, we indicate this by the following diagram, where $\Lm$ consists of just two elements,  $\Lm = \{ 1,2\} $ .
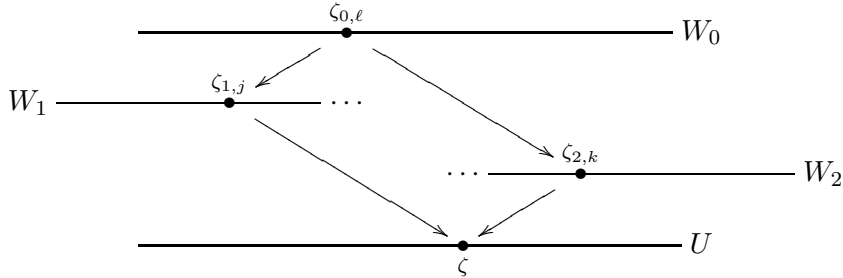
$$
\xymatrix@R=0.3em{   &  \ar@{-}[rrr] & &   \overset{ \zt_{0,\ell}}{\bullet} \ar@{>}[ld] \ar@{>}[rrdd]  \ar@{-}[rrr] & & & W_0 &
\\W_1 \ar@{-}[rrr] &  &  \overset{\zt_{1,j}}{\bullet} \ar@{>}[rrdd]  &  \cdots & & & & &
\\ & &&   & \cdots \ar@{-}[rrr]  & \overset{\zt_{2,k}}{\bullet} \ar@{>}[ld] & & W_2 \\
& \ar@{-}[rrrrr] & & & \underset{\zt}{\bullet} & &  U
}$$
Here $j$ runs through $J_1$, $k$ runs through $J_2$, and $\ell = \{j,k\}$ runs through $J_0$. The arrows stand for the projections onto $W_1$, $W_2,$ $U$.
The preimages of the $\zt_{1,j}$, $\zt_{2,k}$, $\zt_{0,\ell}$ in $\ipi_U(\zt)$ are the members of the partitions $\Pi_1$, $\Pi_2$ and their least common refinement $\Pi_0$. All tangibles $\zt_{1,j}$  and $\zt_{2,k}$ of $W_1$ and $W_2$ have factorization into irreducibles ``closer to UF'' than $\zt$. The tangibles $\zt_{0,\ell}$ of $W_0$ represent the ``least common improvement'' toward unique factorization, compared to the tangibles $\zt_{1,j}$, $\zt_{2,k}$ of $W_1$ and  $W_2$.

\section{Restricted factorization and tangibly  finite monoids}\label{sec:12}

We continue with the setting developed in \S\ref{sec:10} and   \S\ref{sec:11}. If $\zt \in \tT(U)$ is given, then instead of studying all forms where $\zt$ is written as a product of irreducibles, which means to work in the set  $\ipi_U(\zt)$, we choose a monomial $w \in \tT(\mfA)$ and consider only monomials $z \in \ipi_U(\zt)$ with $z | w$, provided that this set of monomials is not empty.

We look at the first instance, where this can be done in a natural way, by taking tangible fiber contractions over $M$ of the supertropical monoid $A_w$ established in \S\ref{sec:9}. Let
\begin{equation}\label{eq:12.1}
  w = \prod_{\iI} x^{r_i}_i \in \tT(\mfA).
\end{equation}
In \S\ref{sec:9} we defined
\begin{equation}\label{eq:12.2}
  A_w := \mfA/ \bigvee_{\iI} \Eq(x^{r_i+1}_i, c^{r_i +1}).
\end{equation}
We now set
\begin{equation}\label{eq:12.3}
  U_w := \mfA/ \bigvee_{\iI} \Eq(x^{r_i+1}_i, c^{r_i +1}) \vee \tlE_U,
\end{equation}
where $\tlE_U$ is given by \eqref{eq:10.4}.

Note that in \eqref{eq:12.2} and \eqref{eq:12.3} only equalizers $\Eq(z, ez)$ with $z \in \tT(\mfA)$ are divided out, while in \S\ref{sec:10} and \S\ref{sec:11} the  equalizers $\Eq(z_1, z_2)$ with $z_1, z_2 \in \tT(\mfA)$, $ez_1 = ez_2$ played center role.

Both $A_w$ and $U_w$ are regarded as subsets of $\mfA$, but they are not submonoids of $\mfA$. Nevertheless they have the same ghost ideal $M$ as $U$, and in particular the same ghost idempotent $e$,
\begin{equation*}\label{eq:12.4}
eA_w = e U_w = eU = e\mfA = M.
\end{equation*}
As in \S\ref{sec:10} we identify the irreducible elements of $\mfA$ with those of $U$
\begin{equation*}\label{eq:12.5}
\Irr(\mfA) = \Irr(U) = \tT_c(U) = \{x_i \pipe \iI \}.
\end{equation*}
It follows that
\begin{equation*}\label{eq:12.5}
\tT_c(A_w) = \tT_c(U_w) = \tT_c(U) = \{x_i \pipe \iI _w \},
\end{equation*}
where $I_w = \{ \iI \pipe r_i > 0\}. $

It is then evident that $A_w$ and $U_w$ both have the properties $A$ and $B$ introduced in \S\ref{sec:10}.
Consequently both $A_w$ and $U_w$ are supertropical semirings and their tangible elements are isolated (Lemma \ref{lem:10.3}).

\begin{thm}\label{thm:12.1}
$ $
\begin{enumerate}\ealph
  \item Given a monomial $w \in \tT(\mfA)$, there is a unique surjective monoid homomorphism
  $\pi_w: A_w\twoheadrightarrow U_w$ with $\pi_w(y) = y$ for all $y\in M$ and $\pi_w(z) = \pi_{\tlU} (z)$.

  \item $U_w = (\tlU)_w$.

  \item Both $\tT(A_w)$ and $\tT(U_w)$ are finite sets.

  \item If $w \in \tT(U)$, then $A_w = U_w$. Otherwise $\tT(U_w) \subsetneqq \tT(A_w)$.

  \item For every $u \in \tT(U_w)$ there exists a unique $z \in \tT(A_w)$ with $\pi_w (z) = u$, and so the  supertropical monoid $U_w$ over $M$ has unique factorization (Remark \ref{rem:10.5}).
\end{enumerate}

\end{thm}
\begin{proof}
(a) is clear form definitions \eqref{eq:12.2}  and \eqref{eq:12.3} of $A_w$ and $U_w$ as quotients of $\mfA$ by the fiberwise equalizers occurring there.
In (b), if we replace $U$ by its cover $\tlU$ in $\mfA$, then the surjective homomorphism $\pi_w$ remains the same map as before, whence $U_w = (\tlU)_w$.
Assertion (c) was dealt in Lemma \ref{lem:10.3}.
To verify (d) and (e) we may replace $U$ by $\tlU$, as consequence of (b). Then they become obvious by (a).
\end{proof}

\begin{cor}\label{cor:12.2}
Assume that $w, w' \in \tT(\mfA)$ and $w' | w$.
\begin{enumerate} \ealph
  \item Then $\pi_{w'}(z) = \pi_{w}(z)$ for every $z \in A_{w'}$.
  \item  There is a unique surjective monoid homomorphism $r_{w,w'}: U_w \to U_{w'}$ such that the diagram
$$
\xymatrix{
A_w  \ar@{>>}[rr]^{\pi_w} \ar@{>>}[d]_{p_{w,w'}} & & U_w \ar@{>>}[d]^{r_{w,w'}}  \\
A_{w'}  \ar@{>>}[rr]^{\pi_{w'}} & & U_{w'}  \\
}
$$
commutes, where $p_{w,w'}$ is the projection \eqref{eq:9.12}. \{The roles of $w$ and $w'$ have been interchanged.\}
\end{enumerate}
\end{cor}
\begin{proof} (a): Clearly $\pi_{w'}(y) = \pi_{w}(y)$ for $y \in M$. Let $z \in \tT(A_{w'}),$ then $z \in \tT(A_w)$. It follows by Theorem \ref{thm:12.1}.(a) that $\pi_{w}(z) =  \pi_{\tlU}(z) = \pi_{w'}(z)$.
\pSkip
(b): A direct consequence of (a) and Theorem \ref{thm:12.1}.(a).
\end{proof}

We are ready for the central result of this section.
\begin{thm}\label{thm:12.4}
Assume that $W$ is a supertropical monoid over $M$ with properties A and~ B (cf. \S\ref{sec:10}), and that
$\tT_e(W) = \{ 1 \}$, $\tT_c(W) \subset \tT_c(U)$. The following are equivalent.
\begin{enumerate} \eroman
  \item $\tT(W)$ is a finite set.
  \item  There exists a monomial $w \in \tT(\mfA)$ such that $W$ is a tangible fiber contraction of $U_w$ over $M$.
\end{enumerate}
\end{thm}
\begin{proof}

(ii)$\Rightarrow$(i): Since  $W$ is a tangible fiber contraction of $U_w$, there is a monomial homomorphism $U_w \to W$ which maps $\tT(U_w)$ onto $\tT(W)$. $\tT(U_w)$ is a finite set, by Theorem ~\ref{thm:12.1}.(c), and thus $\tT(W)$ is finite.

\pSkip
(i)$\Rightarrow$(ii): For every $\iI$ there is a first number $r_i \geq 0$ such that $x_i^{r_1 +1} \in eW$, since $\tT(W)$ is finite, and $r_i > 0 $ only for finitely many $\iI$. Define
$w : = \prod_{\iI} x_i^{r_i} \in \tT(\mfA) $ as in~ \eqref{eq:12.1}. Then $\Irr(W)= \{x_i \pipe r_i > 0 \}$.
For given $\gm \in \tT(W)$ the possible factorizations of $\gm$ into irreducibles of $W$ are given by the monomials  $z \in \tT(\mfA)$ with   $\pi_W(z) = \gm$, $z | w$, and $\pi_U(z) \notin M$, since the $z$ for which $\pi_U(z) \in M$ have been pushed to $M$ in $W$, cf. \eqref{eq:12.3}. Let~ $S$ denotes the set of these monomials
\begin{equation*}\label{eq:12.7}
  \begin{array}{ll}
    S & = \{ z \in \tT(U_w) \pipe \pi_W(z) \in \tT(W), \pi_U(z) \notin M  \}  \\[1mm]
     &  = \{ z \in \tT(\mfA) \pipe \pi_W(z) \in \tT(W), \pi_U(z) \notin M  \}.
  \end{array}
  \end{equation*}
 Hence,  $W$ is the quotient of $U_w$ by the fiberwise equalizer
  \begin{equation*}\label{eq:12.8}
    W = U_w/\Feq(S)
  \end{equation*}
  of the set $S$.
\end{proof}

Together, Theorems \ref{thm:12.1} and \ref{thm:12.4} lead to a classification of all supertropical monoids ~ $W$ over $M$ having properties A and B, which are \textbf{tangibly finite}, i.e., for which the set $\tT(W)$ is finite, together with the embedding of these into $\mfA(I)$ over $M$ for $I$, say, an infinite set. This allows a host of examples of such monoids.

We next pursue the more modest goal to describe, how an explicit list of tangibly finite monoids $W$ over $M$ can be generated (up to monoid isomorphism over $M$), starting with a free monoid  $\mfA = \mfA(I)$ on a finite totally ordered set $\{ x_i \pipe \iI \} $ of commuting variables. These will serve as the  irreducibles of $\mfA$ and $W$. Since for a given monomial $w$ only monomials $z | w$ will play a role, and among these the monomials which are pressed to ghost in $U$, we do not need to specify the supertropical monoid $U$ over $M$ from \S\ref{sec:10}.
It suffices to know which monomials $z|w$ are ghost in every $W$ under consideration. Thus we assume from start that $U = U_w= \tlU_w$.

The monomials $z |w $ with $\pi_U(z) \in M$ may be chosen in $\tT(A_w)$ with the only restriction that, if $\pi_U(z) \in M$, then $\pi_U(zz') \in M$ for any monomial $z'$ such that $zz' |w$.

Instead of describing in general how to obtain a list of all tangibly finite monoids $W$ over $M$, we only give one example. Although it looks simple, it seems to be sufficient to grasp the essentials of a factorization of tangibles into irreducibles in $W$ and covers of~ $W$ in $A_w$.

We take $I = \{1,2 \}$, $\tT_c(U) = \tT_c(U_w) = \{ x_1, x_2\} $, $r_1 = r_2 =3$, hence $w= x_1^3 x_2^3$, and have the following chart of tangibles in $U_w$.
\begin{diagram}\label{diag:12.6} (Tangibles of $U_w$)
$$ \small
\xymatrix@R=0.6em{
& & & x_1^3 \ar@{>}[rd]  \\
& & x_1^2  \ar@{>}[ru] \ar@{>}[rd] & &  x_1^3 x_2  \ar@{>}[rd] \\
& x_1 \ar@{>}[ru] \ar@{>}[rd] & & x_1^2 x_2  \ar@{>}[ru]  \ar@{>}[rd] & & x_1^3 x_2^2 \\
1 \ar@{>}[ru] \ar@{>}[rd] & & x_1 x_2 \ar@{>}[ru] \ar@{>}[rd] & & x_1^2 x_2^2  \ar@{>}[ru]  \ar@{>}[rd]  \\
& x_2 \ar@{>}[rd] \ar@{>}[ru] & & x_1 x^2_2  \ar@{>}[ru]  \ar@{>}[rd] & & x_1^2 x_2^3  \\
& & x_2^2 \ar@{>}[rd]  \ar@{>}[ru]  & & x_1 x_2^3  \ar@{>}[ru] \\
& & &  x_2^3 \ar@{>}[ru] \ . &
}
$$
\end{diagram}

From this chart  we obtain any tangible contraction $W$ if $U_w$ by identifying some pairs of monomials $z_1,z_2$ in the same column by a tangible equalizer $\Eq(z_1,z_2)$  as described in ~ \eqref{eq:10.8}.  We indicate this by a vertical edge between $z_1$ and $z_2$. We then  also have to indicate $(z_1 u, z_2u)$ for $u = x_1$ and $u = x_2$ by an edge. But we do not mark pairs $(z_1,z_3)$ in the same column which are already necessarily identified by  binary equalizers.  For example, if $(z_1,z_2)$, $(z_2,z_3)$ are  marked by an edge, then it is not necessary to mark $(z_1,z_3)$. We write down one of these tables, exhibiting a specific tangibly  finite supertropical monoid ~$W$ over~ $M$.

\begin{diagram}\label{diag:12.7}
$$ \small
\xymatrix@R=0.6em{
& & & x_1^3 \ar@{>}[rd] \ar@{-}[dd] \\
& & x_1^2  \ar@{>}[ru] \ar@{-}[dd] \ar@{>}[rd] & &  x_1^3 x_2 \ar@{-}[dd] \ar@{>}[rd] \\
& x_1 \ar@{>}[ru] \ar@{>}[rd] & & x_1^2 x_2  \ar@{>}[ru]  \ar@{>}[rd] & & x_1^3 x_2^2 \\
1 \ar@{>}[ru] \ar@{>}[rd] & & x_1 x_2 \ar@{>}[ru] \ar@{>}[rd] & & x_1^2 x_2^2 \ar@{-}[dd] \ar@{>}[ru]  \ar@{>}[rd]  \\
& x_2 \ar@{>}[rd] \ar@{>}[ru] & & x_1 x^2_2 \ar@{-}[dd]  \ar@{>}[ru]  \ar@{>}[rd] & & x_1^2 x_2^3  \\
& & x_2^2 \ar@{>}[rd]  \ar@{>}[ru]  & & x_1 x_2^3  \ar@{>}[ru] \\
& & &  x_2^3 \ar@{>}[ru] \ .   &
}
$$
\end{diagram}

From Diagram \ref{diag:12.7} we may read off easily the list of all equalizers $\Eq(z_1, z_2)$ which identify two monomials in $U_w$ (Diagram \ref{diag:12.6}) in $W$, and moreover, which cover $W'$ of $W$ in $U_w$ appears, if we erase one vertical edge in Diagram   \ref{diag:12.7}. For example, if we erase the edge between $x_1^2$ and $x_1 x_2$, no other  equalizers $\Eq(z_1, z_2)$ disappears from the list, i.e., no other vertical edge needs to be erased. But if we erase the edge connecting $x_1^3$ to $x_1^2 x_2$, then the edge between $x_1^2$ and $x_1 x_2$ needs also to be erased, in order to describe $W'$.

When we erase the two most right vertical edges,
then all edges in Diagram   \ref{diag:12.7} need to be erased, and so $W' = U_w$.


\section{The behavior of total and partial splitting covers under transmissions}\label{sec:13}

Assume that $\al: U \onto V$ is a surjective transmission from a supertropical monoid $U$ to a supertropical monoid
$V$ where  $M = eU$ and $N = eV$ are minimal discrete cancellative bipotent semirings. More specifically
$$ M = \{ e, c, c^2, \dots \} \ds{ \dot \cup} \00, \quad c > 0, $$
with all $c^i$ different,
$$ N = \{ e, d, d^2, \dots \} \ds{\dot \cup} \00, \quad d > 0, $$
with all $d^i$ different.
Then $U$ and $V$ have the sets of irreducible elements
\begin{align*}
  \Irr(U) & = \tT_c(U) = \{ x_i \pipe i \in I \},\\
  \Irr(V) & = \tT_d(V) = \{ y_j \pipe j \in J \},
\end{align*}
for some index sets $I$ and $J$ (properties A and B from \S\ref{sec:10}).

As detailed in \S\ref{sec:10}, $U$ and $V$ are fiber contractions of $\mfA(I)$ and $\mfA(J)$ over $M$ and $N$ respectively, i.e., there are multiplicative equivalence relations  $E_U$ and $E_V$ on $\mfA(I)$ and ~ $\mfA(J)$ with associated projections
$$ \pi_U: \mfA(I) \onto U, \qquad \pi_V:\mfA(J) \onto V$$
which restrict to $\id_M$ and~ $\id_N$. We identify the irreducibles $x_i$ and $y_j$ of $\mfA(I)$ and ~ $\mfA(J)$ with irreducibles $x_i$ and $y_j$ in $U$ and $V$, a legal abuse,  since
\begin{equation}\label{eq:13.1}
  [x_i]_{E_U} = \{ x_i \}, \qquad [y_j]_{E_V} = \{y_j\} .
\end{equation}
The monoid $U$ is generated by $M$
 and the $x_i$'s, $i \in I$, while $V$ is generated by $N$
  the $y_j$'s, $j \in J$.

  We introduce the map
  $$\xymatrix@R=3.5em@C=2em{    \tlal := \al \circ \pi_U : \mfA(I)   \ar@{>}[r]^{\quad \qquad \pi_U} & U  \ar@{>}[r]^{\al} &V .} $$
  It is a surjective transmission, since both $\al$ and $\pi_U$ have this property.

  \begin{prop}\label{prop:13.1}
    For any tangible element
    $$ z = \prod_{\iI} x_i^{m_i} \in \mfA(I)$$
    we have
    $$ \tlal(z) = \prod_{\iI}\al( x_i)^{m_i} \in V,$$
    employing  the convention $\pi_U(x_i) = [x_i]_{E_U} = x_i \in U $.
    Either $\tlal(z) \in \tT(V)$ or
    $$ \tlal (z) = e \al(z) = \gm(ez)$$
    where  $\gm:= \al^\nu: M\to N$.
  \end{prop}
  \begin{proof}
    $ \pi_U(z) = \prod_{\iI} \pi_U(x_i)^{m_i} = \prod_{\iI} x_i^{m_i}.$ Applying $\al$ gives
    $\tlal(z)= \prod_{\iI} \al(x_i)^{m_i} \in V$. If this element is not in $\tT(U)$, then
     $$ \tlal(z) = e \prod_{\iI} \al(x_i)^{m_i} = \prod_{\iI} \gm(x_i)^{m_i} = \gm(ez).$$ \vskip -5mm
      \end{proof}
 Assume now that $W$ is a cover of $V = \mfA(J)/E_V$ in $\mfA(J)$, whence $W = \mfA(J)/E_W$ with $E_W \subset E_V$, $V \subset W \subset \mfA(J)$. Assume also  that $H$ is a subset of $\tT(V) $ with
 $\ipi_V(H) \subset W$. Thus $W$ is a total splitting cover of $V$ for $H$. Let
 \begin{equation}\label{eq:13.2}
   \xymatrix@R=2.5em@C=2.8em{   L  \ar@{>}[d]_{\pi_{L,U}} \ar@{>}[r]^{\bt} & W  \ar@{>}[d]^{\pi_{W,V}}\\
   U \ar@{>}[r]^{\al} &V }
 \end{equation}
 denote the fiber contraction of $U \overset{\al}{\To} V$ and $W \overset{\al}{\To} V$ over $V$, i.e., the pulback of $\al$ and~ $\pi_{W,V}$ in the category of supertropical monoids with transmissions. \{In particular $\bt$ is again a surjective transmissions.\}

\begin{thm}\label{thm:13.2}Let  $L$ be a total splitting cover of $U$ for the set $\ipi(H) \subset \tT(U)$. If the total splitting cover $W$ of $V$ for $H$ is minimal, then the splitting cover $L$ of $U$ for
$\ipi(H)$ is also minimal. The ghost maps $\al^\nu$ and $\bt^\nu$ are the same homomorphism $\gm : M \to N$.

\end{thm}
\begin{proof}
  $L$ may be viewed as the set of pairs $(z_1, z_2)$ with $z_1 \in U$, $z_2 \in W$ and $\al(z_1) = \pi_{W,V} (z_2)$. Considering this and recalling the relevant setup from \S\ref{sec:9},  the claims can be verified in a straightforward way.
\end{proof}

\end{document}